\magnification=1150
\overfullrule=0pt
\parindent=0pt
\pageno=1

\font\rmsec=cmbx10 scaled 1300

\input epsf.def
\input epsf.sty 
\input amssym.def
\input amssym.tex
\input xy
\xyoption{all}

\def\address{
I{\sevenrm NSTITUTE OF} M{\sevenrm ATHEMATICS}

T{\sevenrm HE} H{\sevenrm EBREW} U{\sevenrm NIVERSITY}

Givat Ram Campus 

91904 Jerusalem -- Israel

\pec 

E-mail: {\tt remy@math.huji.ac.il}

Homepage: http:$/ \! \! /$www.ma.huji.ac.il/$^\sim$remy}

\font \tengoth=eufm10 at 12pt
\font \sevengoth=eufm7
\font \fivegoth=eufm5
\newfam\gothfam
\textfont \gothfam=\tengoth
\scriptfont \gothfam=\sevengoth
\scriptscriptfont \gothfam=\fivegoth
\def\goth{\fam\gothfam\tengoth}

\def\pec{\vskip2mmplus1mmminus1mm}
\def\gec{\vskip4mmplus1mmminus1mm}
\def\pech{\hskip3mmplus1mmminus1mm}

\def\og{\leavevmode\raise.3ex\hbox{$\scriptscriptstyle\langle\!\langle\,$}}
\def\fg{\leavevmode\raise.3ex\hbox{$\scriptscriptstyle\,\rangle\!\rangle\ $}}
\def\trait{--- \pech }
\def\tch#1{#1\mkern2.5mu\check{}}
\def\cqfd{\nobreak\hfill $\square$ \goodbreak}


\def\iso{\buildrel\sim \over \rightarrow}


\def\N{\hbox {\bf N}}
\def\C{\hbox {\bf C}}
\def\K{\hbox {\bf K}}
\def\R{\hbox {\bf R}}
\def\Z{\hbox {\bf Z}}
\def\E{\hbox {\bf E}}
\def\fq{\hbox {\bf F}_q}
\def\fqd{\hbox {\bf F}_{q^2}}
\def\pz{\hbox {\sevenbf Z}}
\def\pe{\hbox {\sevenbf E}}
\def\pn{\hbox {\sevenbf N}}
\def\pr{\hbox {\sevenbf R}}
\def\pk{\hbox {\sevenbf K}}
\def\ks{\hbox {${\bf K}_s$}}

\def\kb{\overline {\hbox{\bf K}}}
\def\pkb{\overline {\hbox{\sevenbf K}}}


\def\gd{\g_{\cal D}}
\def\ugd{{\cal U}\g_{\cal D}}
\def\ud{{\cal U}_{\cal D}}
\def\uo{{\cal U}_0}
\def\st{\hbox{\rm St}_A}
\def\gkm{\underline G_{\cal D}}
\def\uab{{\cal U}_{[ a;b ]_{\hbox{\sevenrm lin}}}}
\def\autf{\aut_{filt}({\cal U}_{\cal D})}
\def\ug{\underline G}
\def\uk{{\cal U}_{\pk}}

\def\ukb{\bigl( \ud \bigr)_{\pkb}}


\def\ip{{\cal I}_+}
\def\im{{\cal I}_-}

\def\A{\hbox {\bf A}}
\def\ak{\A_{\pk}}
\def\ct{\overline {\cal C}}
\def\ict{{\cal C}}
\def\fr{F^\natural} 
\def\lr{L^\natural} 
\def\arel{a^\natural} 
\def\br{b^\natural} 


\def\g{\hbox {\goth g}}
\def\ad{\hbox {\rm ad}}
\def\Ad{\hbox {\rm Ad}}


\def\aut{\hbox {\rm Aut}}
\def\paut{\hbox {\sevenrm Aut}}
\def\Hom{\hbox {\rm Hom}}


\def\trdl{\hbox{\rm (TRD)}_{\hbox {\sevenrm lin}}}

\centerline{\rmsec KAC-MOODY GROUPS:}

\gec

\centerline{\rmsec SPLIT AND RELATIVE THEORIES. LATTICES}

\gec

\centerline{\bf B{\sevenbf ERTRAND} R\'EMY}

\vskip 10mm

A{\sevenrm BSTRACT}.---~{\it In this survey article, we recall some facts about split
Kac-Moody groups as defined by J. Tits, describe their main properties and then propose an
analogue of Borel-Tits theory for a non-split version of them.  The main result is a Galois
descent theorem, i.e., the persistence of a nice combinatorial structure after
passing to rational points. We are also interested in the geometric point of view, namely
the production of new buildings admitting $($nonuniform$)$ lattices}. 

\vskip 10mm

\centerline{CONTENTS}

\gec

INTRODUCTION

\pec

1. ABSTRACT GROUP COMBINATORICS AND TWIN BUILDINGS

\pec

1.1. Root system and realizations of a Coxeter group 

1.2. Axioms for group combinatorics 

1.3. Moufang twin buildings

1.4. A deeper use of geometry

1.5. Levi decompositions

\pec

2. SPLIT KAC-MOODY GROUPS

\pec

2.1. Tits functors and  Kac-Moody groups

2.2. Combinatorics of a Kac-Moody group

2.3. The adjoint representation

2.4. Algebraic subgroups

2.5. Cartan subgroups

\pec

3. RELATIVE KAC-MOODY THEORY

\pec

3.1. Definition of forms

3.2. The Galois descent theorem 

3.3. A Lang type theorem

3.4.  Relative apartments and relative links

3.5. A classical twin tree

\pec

4. HYPERBOLIC EXAMPLES. LATTICES

\pec

4.1. Hyperbolic examples 

4.2. Analogy with arithmetic lattices

\vfill

{\it $1991$ Mathematics Subject Classification:} 20E42, 51E24, 20G15, 22E65. 


\eject

\centerline{\bf INTRODUCTION} 

\gec

{\it Historical sketch of Kac-Moody theory.---~} Kac-Moody theory was initiated in
1968, when V. Kac and R. Moody defined independently infinite-dimensional Lie algebras
generalizing complex semisimple Lie algebras. Their  definition is based on Serre's
presentation theorem describing explicitly the latter
(finite-dimensional) Lie algebras [Hu1, 18.3]. A natural
question then is to integrate Kac-Moody Lie algebras as Lie groups integrate real Lie
algebras, but this time in the infinite-dimensional setting.  This difficult problem led
to several propositions. In characteristic 0, a satisfactory
approach consists in seeing them as subgroups in the automorphisms  of the corresponding
Lie algebras [KP1,2,3]. This way, V. Kac and D. Peterson developed the structure theory of
Kac-Moody algebras in complete analogy with the classical theory: intrinsic definition
and conjugacy results for Borel (resp. Cartan) subgroups, root decomposition with abstract
description of the root system... Another aspect of this work is the construction of
generalized Schubert varieties. These algebraic varieties enabled O. Mathieu to get a
complete generalization of the character formula in the Kac-Moody framework [Mat1]. To this
end, O. Mathieu defined Kac-Moody groups over arbitrary fields in the formalism of
ind-schemes [Mat2]. 

\pec

{\it Combinatorial approach.---~} Although the objects above -- Kac-Moody groups
and Schubert varieties -- can be studied in a nice algebro-geometric context, we will
work with groups arising from another, more combinatorial viewpoint. All of this work is due
to J. Tits [T4,5,6,7], who of course contributed also to the previous problems.  The aim is
to get a much richer combinatorial structure for these groups. This led J. Tits to 
the notion of \og Root Group Datum\fg axioms [T7] whose geometric counterpart is the theory of
Moufang twin buildings. The starting point of the construction of Kac-Moody groups [T4] is
a generalization of Steinberg's presentation theorem [Sp, Theorem 9.4.3] which concerns
simply connected semisimple algebraic groups. In this context, the groups ${\rm
SL}_n(\K[t,t^{-1}])$ are Kac-Moody groups obtained via Tits' construction. 

\pec 

{\it Relative Kac-Moody theory in characteristic $0$.---~} So far, the
infinite-dimensional objects alluded to were analogues of split Lie algebras and split
algebraic groups. Still, it is known that by far not all interesting algebraic groups  are
covered by the split theory -- just consider the simplest case of the multiplicative group
of a quaternion skewfield. This is the reason why Borel-Tits theory [BoT]  is so important: it deals with algebraic groups over arbitrary fields $\K$, and the main results are a
combinatorial structure theorem for $\K$-points, conjugacy theorems for minimal
$\K$-parabolic subgroups (resp. for maximal $\K$-split tori). This theory calls for a
generalization in the Kac-Moody setting: this work was achieved in the characteristic 0
case by G. Rousseau [Rou1,2,3; B$_3$R]. 

\pec

{\it Two analogies.---~} The example of the groups ${\rm SL}_n(\K[t,t^{-1}])$ is
actually a good guideline since it provides at the same time another analogy for Kac-Moody
groups. Indeed, over finite fields $\K=\fq$, the former ones are arithmetic groups in the
function field case. Following the first analogy -- a Kac-Moody group is an
infinite-dimensional reductive group, a large part of the present paper describes the
author's thesis [R\'e2] whose basic goal is to define $\K$-forms of Kac-Moody groups with
no assumption on the groundfield $\K$. This prevents from using the viewpoint of
automorphisms of Lie algebras as in the works previously cited, but the analogues of the
main results of Borel-Tits theory are proved. This requires first to reconsider Tits'
construction of Kac-Moody groups and to prove results in this (split) case, interesting in
their own right. On the other hand, the analogy with arithmetic groups will also be
discussed so as to see Kac-Moody groups as discrete groups -- lattices for $\fq$ large
enough -- of their geometries. 

\pec 

{\it Tools.---~} Let us talk about tools now, and start with the main difficulty: no 
algebro-geometric structure is known for split Kac-Moody groups as defined by J. Tits. The
idea is to replace this structure by two well-understood actions. The first one is not so
mysterious since it is linear: it is in fact the natural generalization of the adjoint
representation of algebraic groups. It plays a crucial role because it enables to endow a large family of
subgroups with a structure of algebraic group. 
The second kind of action involves buildings: it is a usual topic in group theory to define 
a suitable geometry out of a given group to study it (use of Cayley graphs, of boundaries... ) 
Kac-Moody groups are concerned by building techniques thanks to their nice combinatorial 
properties, refining that of $BN$-pairs (for an account on the general use of buildings 
in group theory, see [T2]). There exists some kind of (non-univoque) correspondence between
buildings and $BN$-pairs [Ron \S 5], but it is too general to provide precise 
information in specific situations. Some refinements are to be adjusted accordingly 
-- see the example of Euclidean buildings and \og Valuated Root Data\fg in Bruhat-Tits theory [BrT1]. 
As already said, in the Kac-Moody setting the refinements consist in requiring 
the \og Root Group Datum\fg axioms at the group level and in working with twin buildings 
at the geometric level. Roughly speaking, a twin building is the datum of two buildings 
related by opposition relations between chambers and apartments. 

\pec 

{\it Organization of the paper.---~} This article is divided into
four parts. Part 1 deals with group combinatorics in a purely abstract
context. Since the combinatorial axioms will be satisfied by both split and non-split
groups, it is an efficient way to formalize properties shared by them. There is described
the geometry of twin buildings and the corresponding group theoretic axioms. The aim is
to obtain two kinds of Levi decomposition of later interest. Part 2 describes
the split theory of Kac-Moody groups. In particular, we explain why the adjoint
representation can be seen as a substitute for a global algebro-geometric structure. An
illustration of this is a repeatedly used argument combining negative curvature and
algebraic groups arguments. Part 3 presents the relative theory of almost split Kac-Moody
groups. A sketch of the proof of the structure theorem for rational points is given. The
particular case of a finite groundfield is considered, as well as a classical example of a
twisted group leading to a semi-homogeneous twin tree. At last, part 4 adopts the
viewpoint of discrete groups. We first show that Kac-Moody theory enables to produce
hyperbolic buildings (among many other possibilities) and justify why these geometries
are particularly interesting. Then we show that Kac-Moody groups or their spherical
parabolic subgroups over a finite groundfield are often lattices of their buildings. 
This leads to an analogy with arithmetic lattices over function fields.  

\pec 

The assumed knowledge for this article consists of general facts from building theory and
algebraic groups. References for buildings are K. Brown's book [Br2] for the apartment
systems viewpoint and M. Ronan's book [Ron] for the chamber systems one. Concerning
algebraic groups, the recent books  [Bo] and [Sp] are the main references dealing with
relative theory. 

\pec 

{\it Acknowledgements.---~} This work presents a Ph.D. prepared under the
supervision of G. Rousseau. It is a great pleasure to thank him, as well as M.
Bourdon who drew my attention to the discrete groups viewpoint. I am very grateful to H.
Abels, P. Abramenko, Th. M\"uller for their kindness. It was a pleasure and a nice
experience to take part in the conference \og Groups: geometric and combinatorial
aspects\fg organized by H. Helling and the latter, 
and to be welcomed at the SFB 343 (University of Bielefeld). 

\vfill\eject

\centerline{\bf 1. ABSTRACT GROUP COMBINATORICS AND TWIN BUILDINGS}

\gec

The aim of this section is to provide all the abstract background we
will need to study split and almost split Kac-Moody groups. In \S 1.1 is
introduced the root system of a Coxeter group. These roots will index
the group combinatorics of the \og Twin Root Datum\fg axioms presented in \S
1.2. Then we describe in \S 1.3 the geometric side of the (TRD)-groups, 
that is the twin buildings on which they operate. Section 1.4 is
dedicated to the geometric notions and realizations to be used later. At last,
\S 1.5 goes back to group theory providing some semi-direct decomposition
results for distinguished classes of subgroups. 

\gec

\centerline{\bf 1.1. Root system and realizations of a Coxeter group}

\pec

The objects we define here will be used to index group theoretic axioms (1.2.A) and to
describe the geometry of buildings (1.3.B). 

\pec

{\bf A. Coxeter complex. Root systems. } Let $M=[M_{st}]_{s,t \in S}$ denote a Coxeter
matrix with associated Coxeter system $(W,S)$. For our purpose, it is sufficient to suppose
the canonical generating set $S$ finite. The group $W$ admits the following presentation: \pec

\centerline{$W = \langle s \!\in \! S \mid (st)^{M_{st}}=1$ whenever $M_{st} < \infty 
\rangle$.}

\pec

We shall use the {\it length function } $\ell: W \to \N$ defined w.r.t. $S$. 

\pec 

The existence of an abstract simplicial complex acted upon by $W$ is the starting point
of the definition of buildings of type $(W,S)$ in terms of apartment systems. This
complex is called the {\it Coxeter complex } associated to $W$, we will denote it by
$\Sigma(W,S)$ or $\Sigma$. It describes the combinatorial geometry of \og slices\fg in a
building of type $(W,S)$ -- the so-called {\it apartments}.  The abstract complex $\Sigma$
is made of translates of the {\it special subgroups} $W_J:= \langle J \rangle$, $J \subset
S$, ordered by inclusion [Br2, p.58-59]. 

\pec

The {\it root system } of $(W,S)$ is defined by means of the length function $\ell$
[T4, \S 5].  The set $W$ admits a $W$-action via left translations. {\it Roots } are
distinguished halves of this $W$-set, whose elements will be called {\it chambers}. 

\pec 

{\bf Definition.---~} \it 
{\rm (i)} The {\rm simple root } of index $s$ is the half  
$\alpha_s:=\{ w \! \in \! W \mid \ell(sw) > \ell(w) \}$. 

{\rm (ii) } A {\rm root} of $W$ is a half of the form $w\alpha_s$, $w \! \in \! W, 
s \! \in \! S$. The set of roots will be denoted by $\Phi$, it admits an obvious
$W$-action. 

{\rm (iii) } A root is called {\rm positive } if it contains
$1$; otherwise, it is called {\rm negative}. Denote by $\Phi_+$
$($resp. $\Phi_-)$ the set of positive $($resp. negative$)$ roots of $\Phi$. 

{\rm (iv) } The {\rm opposite } of a root is its complement. 
\rm\pec 

{\it Remark.---~} The opposite of the root $w\alpha_s$ is indeed a root since it is
$ws\alpha_s$. 

\pec

The next definitions are used for the group combinatorics presented in 1.2.A. 

\pec 

{\bf Definition.---~} \it 
{\rm (i) } A pair of roots $\{\alpha; \beta\}$ is called {\rm prenilpotent } if both
intersections $\alpha \cap \beta$ and $(-\alpha) \cap (-\beta)$ are
nonempty. 

{\rm (ii) } Given a prenilpotent pair of roots $\{\alpha; \beta\}$,
the {\rm interval }$[\alpha; \beta]$ is by definition the set of roots
$\gamma$ with $\gamma \supset \alpha \cap \beta$ and 
$(-\gamma) \supset (-\alpha) \cap (-\beta)$. We also set 
$]\alpha; \beta[:=[\alpha; \beta] \setminus \{\alpha;\beta\}$. 
\rm\pec 

{\bf B. The Tits cone. Linear refinements. } We introduce now a fundamental
realization of a Coxeter group: the {\it Tits cone}. It was first defined in
[Bbk, V.4] to which we refer for proofs.  This approach will also allow to consider
roots as linear forms [Hu2, II.5].  We keep the Coxeter system $(W,S)$ with Coxeter
matrix $M=[M_{st}]_{s,t \in S}$. Consider the real vector space $V$ over the symbols 
$\{ \alpha_s \}_{s \in S}$ and define the {\it cosine matrix } $A$ of $W$ by 
$A_{st}:= -\cos(\pi/M_{st})$: this is the matrix of a bilinear form
$B$ (w.r.t. the basis $\{ \alpha_s \}_{s \in S}$). To each $s$ of $S$ is associated the
involution $\sigma_s: \lambda \mapsto \lambda - 2B(\alpha_s,\lambda)\alpha_s$, and 
the assignment $s \mapsto \sigma_s$ defines a faithful representation of $W$. 
The {\it (positive) half-space } in $V^*$ of an element $\lambda$ in $V$ is denoted by 
$D(\lambda)$: $D(\lambda):= \{ x \! \in \! V^* \mid \lambda(x)>0 \}$. We denote its
boundary -- the kernel of $\lambda$ -- by $\partial\lambda$. 

\pec 

{\bf Definition.---~} \it 
{\rm (i) } The {\rm standard chamber } $c$ is the simplicial cone $\displaystyle 
\bigcap_{s \in S} D(\alpha_s)$ of elements of $V^*$ on which all linear forms $\alpha_s$
are positive. A {\rm chamber } is a $W$-translate of $c$. 

{\rm (ii) } The {\rm standard facet of type $J$}, $J \subset S$, is 
$\displaystyle F_I:= \bigcap_{s \in I} \partial\alpha_s \cap 
\bigcap_{s \in S \setminus I} D(\alpha_s)$. A {\rm facet of type $J$ } is a 
$W$-translate of $F_J$. 

{\rm (ii) } The {\rm Tits cone $\ct$ } of $W$ is the union of the closures of all chambers
$\displaystyle \bigcup_{w \in W} w \overline c = \bigsqcup_{w \in W, J \subset S} wF_J$. 
\rm\pec

A study of the action of $W$ on $\ct$ shows that the type is well-defined. A facet is {\it
spherical } if it is of type $J$ with $W_J$ finite. Facets of all types are
represented here, as simplicial cones. The simplicial complex so obtained is not locally
finite in general, but its interior $\ict$ is. In fact, a facet is of nonspherical type if
and only if it lies in the boundary of the Tits cone. Further properties of the cone $\ct$
are available in [V]. 

\pec 

{\it Example.---~} In the case of an affine reflexion group, the Tits cone is made of
the union of an open half-space and the origin, which is the only nonspherical facet. The
affine space in which the standard representation of the group is defined is just the
affinisation of this cone [Hu2, II.6]. 

\pec

The viewpoint of linear forms for roots enables to introduce another -- more restrictive --
notion of interval of roots. 

\pec 

{\bf Definition.---~} \it 
Given two roots $\alpha$ and $\beta$, the {\rm linear interval } they
define is the set $[\alpha;\beta]_{\hbox{\sevenrm lin}}$ of positive linear
combinations of them, seen as linear forms on $V$. 
\rm\pec

{\it Figure.---~} 

\centerline{\epsfysize=30mm $$\epsfbox{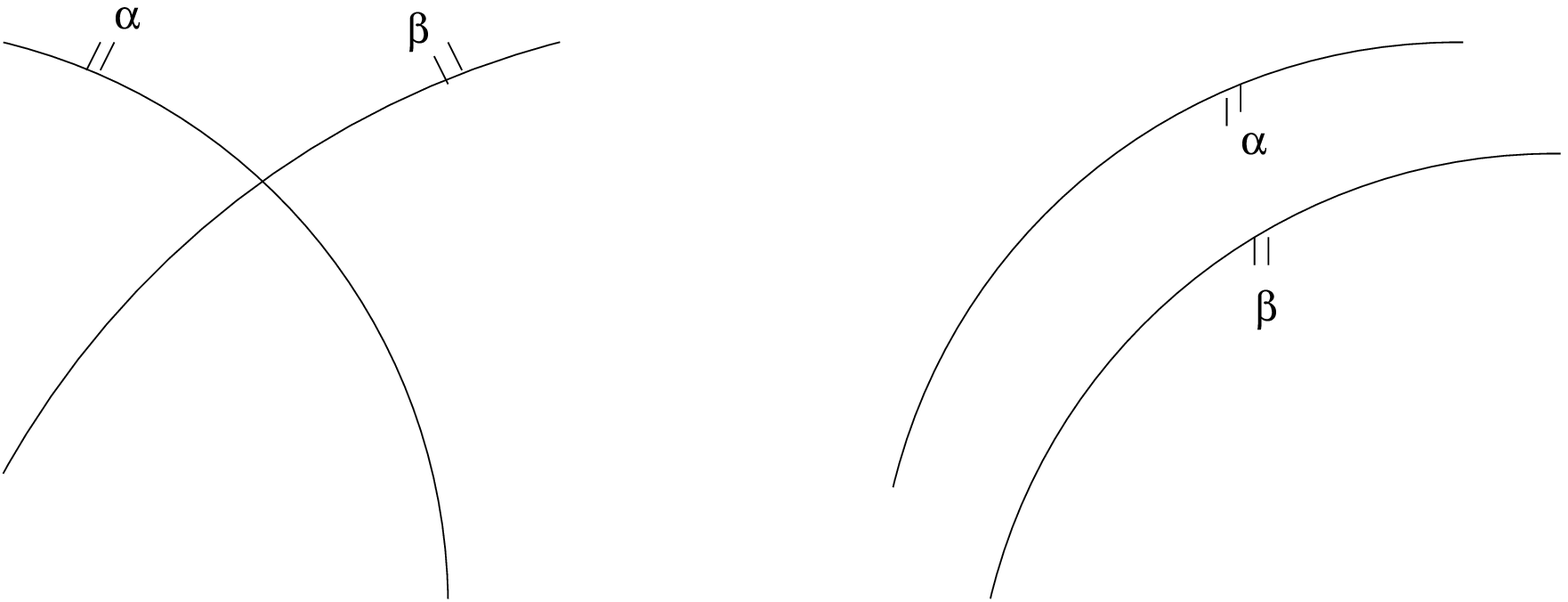}$$}

The picture above illustrates a general fact: there are two ways to be
prenilpotent. Either the walls of the roots intersect along a
spherical facet -- and the four pairs of roots of the form 
$\{ \pm\alpha; \pm\beta \}$ are prenilpotent, or a root contains the other 
-- and only two pairs among the four are prenilpotent. 

\pec

{\it Remark.---~} For a prenilpotent pair of roots, one has
$[\alpha;\beta]_{\hbox{\sevenrm lin}} \subset [\alpha;\beta]$, still the notions do not
coincide, even in the Kac-Moody situation. G. Rousseau indicated an example of strict
inclusion provided by a tesselation of the hyperbolic plane ${\Bbb H}^2$ [R\'e2, 5.4.2]. 

\pec

{\bf C. Geometric realizations. } We will try to use geometry as much as possible
instead of abstract set-theoretic structures. As an example, we will often 
represent Coxeter complexes by {\it polyhedral complexes } [BH]. We are interested in
such spaces with the following additional properties. 

\pec

(i) The complex is labelled by a fixed set of subsets of $S$ -- the {\it types}. We call
{\it facets } the polyhedra inside.  Codimension $1$ facets are called {\it panels},
maximal facets are called {\it chambers}. 

(ii) There is a countable family of codimension $1$ subcomplexes 
-- the {\it walls } -- w.r.t. which are defined involutions -- the {\it reflections}. 
A reflection fixes its wall and stabilizes the whole family of them. 

(iii) For a chamber $c$, there is a bijection between the set of generators $S$ and the
walls supporting the panels of $c$. The corresponding reflections define a faithful
representation of $W$ by label-preserving automorphisms of the complex. 

(iv) The $W$-action is simply transitive on chambers. 

\pec

{\it Remarks.---~} 1. The Tits cone satisfies all the conditions above, with
simplicial cones as facets instead of polyhedra. 

2. We may use geometric realizations where only spherical types appear. This is the case in
the examples below. 

\pec

{\it Examples.---~} 1. The simplest example with infinite Coxeter
group is given by the real line and its tesselation by the segments
defined by consecutive integers. It is acted upon by the infinite
dihedral group $D_\infty$, and the corresponding buildings are trees. 

2. Another famous example associated to an affine reflection
group comes from the tesselation of the Euclidean plane by
equilateral triangles. Buildings with this geometry as apartments are called
$\widetilde A_2$- or {\it triangle buildings}. They are interesting
because, even if they belong to the well-known class of {\it Euclidean
buildings}, many of them do not come from Bruhat-Tits theory.

3. A well-known way to construct concretely a Coxeter group with a realization of its
Coxeter complex is to apply Poincar\'e polyhedron theorem [Mas, IV.H.11]. This 
works in the framework of spherical, Euclidean or hyperbolic geometry; we just have to
consider reflections w.r.t. to a suitable polyhedron. For a tiling of a hyperbolic space
${\Bbb H}^n$, the corresponding buildings are called {\it hyperbolic}, {\it Fuchsian } in
the two-dimensional case.

\pec 

{\it Figure.---~} 

\centerline{\epsfysize=40mm $$\epsfbox{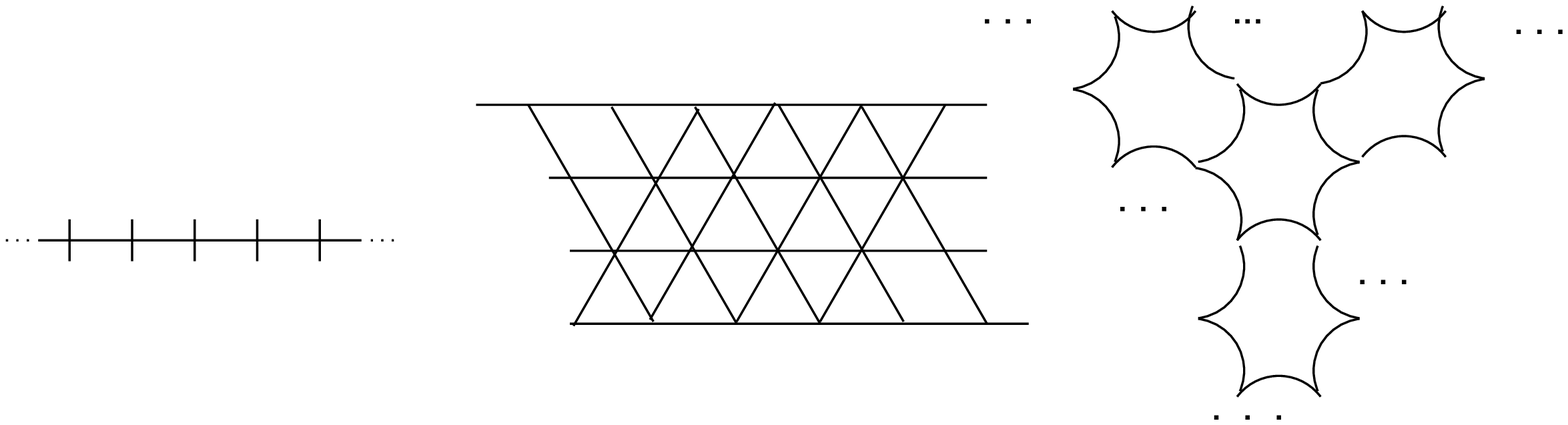}$$}

\centerline{\bf 1.2. Axioms for group combinatorics}

\pec

We can now give axioms refining $BN$-pairs and adapted to the Kac-Moody situation. 

\pec

{\bf A. The (TRD) axioms. } The axioms listed below are indexed by the set
of roots $\Phi$ of the Coxeter system $(W,S)$. It is a slight modification of axioms
proposed in [T7]. That it implies the group combinatorics introduced by
V. Kac and D. Peterson [KP2] follows from [R\'e2, Th\'eor\`eme 1.5.4], which
elaborates on [Ch] and [T3]. 

\pec 

{\bf Definition.---~} \it 
Let $G$ be an abstract group containing a subgroup $H$. 
Suppose $G$ is endowed with a family $\{ U_\alpha \}_{\alpha \in \Phi}$ of subgroups 
indexed by the set of roots $\Phi$, and define the subgroups 
$U_+:=\langle U_\alpha \mid \alpha \! \in \! \Phi_+ \rangle$ and 
$U_-:=\langle U_\alpha \mid \alpha \! \in \! \Phi_- \rangle$. Then, the triple 
$\bigl( G, \{ U_\alpha \}_{\alpha \in \Phi}, H \bigr)$ 
is said to satisfy the {\rm (TRD) axioms } if the following conditions are satisfied.

\pec 

\pech {\rm (TRD0) } Each $U_\alpha$ is nontrivial and normalized by $H$. 

\pech {\rm (TRD1) } For each prenilpotent pair of roots 
$\{ \alpha;\beta \}$, the commutator subgroup $[U_\alpha,U_\beta]$ is
contained in the subgroup $U_{]\alpha,\beta[}$ generated by the $U_\gamma$'s, with 
$\gamma \! \in  ]\alpha,\beta[$. 

\pech {\rm (TRD2) } For each $s$ in $S$ and $u$ in
$U_{\alpha_s} \setminus \{ 1 \}$,
there exist uniquely defined $u'$ and $u''$ in $U_{-\alpha_s} \setminus \{ 1 \}$ such that
$m(u):=u'uu''$ conjugates $U_\beta$ onto $U_{s\beta}$ for every root
$\beta$. Besides, it is required that for all $u$ and $v$ 
in $U_{\alpha_s} \setminus \{ 1 \}$, one should have $m(u)H=m(v)H$. 

\pech {\rm (TRD3) } For each $s$ in $S$, 
$U_{\alpha_s} \not\subset U_-$ and $U_{-\alpha_s} \not\subset U_+$. 

\pech {\rm (TRD4) } $G$ is generated by $H$ and the $U_\alpha$'s. 

\pec 

Such a group will be referred to as a {\rm (TRD)-group}. It will
be called a {\rm $\trdl$-group } or said to satisfy the {\rm $\trdl$ axioms }
if {\rm (TRD1) } is still true after replacing intervals by linear ones.
\rm\pec

{\it Remarks.---~} 1. A consequence of Borel-Tits theory [BoT, Bo, Sp] is that
isotropic reductive algebraic groups satisfy $\trdl$ axioms. We will see in 2.3.A
that so do split Kac-Moody groups. The case of nonsplit Kac-Moody groups is
the object of the Galois descent theorem  -- see 3.2.A. 

2. The case of algebraic groups suggests to take into account more carefully
proportionality relations between roots seen as linear forms. It is indeed possible
to formalize the difference between reduced and non-reduced infinite root systems, and to
derive refined $\trdl$ axioms [R\'e2, 6.2.5]. 

\pec

{\bf B. Main consequences. } We can derive a first list of properties for
a (TRD)-group $G$.

\pec

{\it Two $BN$-pairs.---~} The main point is the existence of two $BN$-pairs in the
group $G$. Define the {\it standard Borel subgroup of sign $\epsilon$ }
to be $B_\epsilon:=HU_\epsilon$. The subgroup $N < G$ is by definition generated by $H$ and
the $m(u)$'s of axiom (TRD2). Then, one has 

\pec 

\centerline{$\displaystyle H = \bigcap_{\alpha \in \Phi} N_G(U_\alpha) = B_+ \cap N = B_-
\cap N$,}

\pec 

and $(G,B_+,N,S)$ and $(G,B_-,N,S)$ are $BN$-pairs sharing
the same Weyl group $W=N/H$. As $B_+$ and $B_-$
are not conjugate, the positive and the negative $BN$-pairs do not carry the same
information. A conjugate of $B_+$ (resp. $B_-$) will be called a {\it
positive (resp. negative) Borel subgroup}. 

\pec

{\it Refined Bruhat and Birkhoff decompositions.---~} A formal consequence of the
existence of a $BN$-pair is a Bruhat decomposition for the
group. In our setting, the decomposition for each sign can be made more
precise. For each $w \! \in \! W$, define the subgroups
$U_w:= U_+ \cap wU_-w^{-1}$ and $U_{-w}:= U_- \cap wU_+w^{-1}$. The {\it refined Bruhat
decompositions } are then [KP2, Proposition 3.2]: $$G = \bigsqcup_{w \in W} U_wwB_+ \quad \hbox{\rm and } 
\quad G = \bigsqcup_{w \in W} U_{-w}wB_-,$$ 

\pec 

with uniqueness of the first factor. A third decomposition involves both
signs and will be used to define the twinned structures (1.3.A). More precisely,
the {\it refined Birkhoff decompositions } are [KP2, Proposition 3.3]: \pec 

$$G = \bigsqcup_{w \in W} (U_+ \cap wU_+w^{-1})wB_-
= \bigsqcup_{w \in W} (U_- \cap wU_-w^{-1})wB_+,$$

\pec 

once again with uniqueness of the first factors. 

\pec

{\it Other unique writings.---~} 
Another kind of unique writing result is valid for the groups $U_{\pm w}$. 
For each $z \! \in \! W$, define the (finite) sets of roots 
$\Phi_z:= \Phi_+ \cap z^{-1}\Phi_-$ and $\Phi_{-z}:= \Phi_- \cap z^{-1}\Phi_+$. 
Then, the group $U_w$ (resp. $U_{-w}$) is in bijection with the set-theoretic
product of the root groups indexed by $\Phi_{w^{-1}}$ (resp. $\Phi_{-w^{-1}}$) for a suitable 
(cyclic) ordering on the latter set [T4, proposition 3 (ii)], [R\'e2, 1.5.2].  

\pec 

{\it Two buildings.---~} Let us describe now how to construct a building
out of each $BN$-pair, the twin structure relating them being
defined in 1.3.A. Fix a sign $\epsilon$ and consider the
corresponding $BN$-pair $(G,B_\epsilon,N,S)$. 
Let $d_\epsilon: G/B_\epsilon \times G/B_\epsilon \rightarrow W$ be the application which
associates to $(gB_\epsilon,hB_\epsilon)$ the element $w$ such that $B_\epsilon g^{-1}h
B_\epsilon=B_\epsilon w B_\epsilon$. Then, $d_\epsilon$ is a {\it $W$-distance } making 
$(G/B_\epsilon,d_\epsilon)$ a building [Ron, 5.3]. The {\it standard apartment
of sign $\epsilon$ } is $\{ wB_\epsilon \}_{w \in W}$, the relevant apartment system
being the set of its $G$-transforms. A {\it facet } of type $J$ is a translate
$gP_{\epsilon,J}$ of a {\it standard parabolic subgroup } $P_{\epsilon,J}:= B_\epsilon W_J
B_\epsilon$. 

\pec 

{\bf C. } {\it Examples.---~} The most familiar examples of groups enjoying the
properties above are provided by Chevalley groups over Laurent polynomials. 
These groups are Kac-Moody groups of affine type. We briefly describe the case of the
special linear groups ${\rm SL}_n(\K[t,t^{-1}]), n \geq 2$. From the Kac-Moody
viewpoint, the groundfield is $\K$. The buildings involved are Bruhat-Tits. They are
associated to the
$p$-adic Lie groups ${\rm SL}_n \bigl( \fq(\!(t)\!) \bigr)$ and ${\rm SL}_n \bigl(
\fq(\!(t^{-1})\!) \bigr)$ respectively in the case of a finite groundfield $\fq$. 
The Weyl group is the affine reflection group ${\cal S}_{n} \ltimes \Z^{n-1}$. 
The Borel subgroups are

$$B_+:= \{ M \in \pmatrix{\xymatrix{\K[t] \ar@{.}[dr] & \K[t] \\ t\K[t] & \K[t]}} 
\mid {\hbox {\rm det} } M=1 \};$$

$$B_-:= \{ M \in \pmatrix{\xymatrix{ \K[t^{-1}] \ar@{.}[dr]&
t^{-1}\K[t^{-1}] 
\\ \K[t^{-1}] & \K[t^{-1}]}}
\mid {\hbox {\rm det} } M=1 \}.$$

As subgroup $H$, we take the standard Cartan subgroup $T$ of 
${\rm SL}_n (\K)$ made of diagonal matrices with coefficients in
$\K^\times$ and determinant 1. For $1 \leq i \leq n-1$, the
monomial matrices with $\pmatrix{0&1\cr-1&0\cr}$ in position $(i,i+1)$ on the diagonal 
(and 1's everywhere else on it) lift the $n-1$ simple reflexions 
generating the finite Weyl group of ${\rm SL}_n (\K)$. The last reflexion \og responsible
for the affinisation\fg is lifted by 

\pec

\centerline{$N_n:=\pmatrix{0&0&-t^{-1} \cr 
0&\pmatrix{\xymatrix{1\ar@{.}[dr]&0\\0&1}}&0 \cr
t&0&0 \cr}$.} 

\pec

The situation is the same for root groups: besides the simple root groups of 
${\rm SL}_n (\K)$, one has to add 

\centerline{$U_n:=\{ u_n(k):= \pmatrix{\xymatrix{
1\ar@{.}[dr]&0&0 \\ 0&1\ar@{.}[dr]&0 \\ kt&0&1}} \mid k\!\in\!\K \}$,}

\pec 

to get the complete family of subgroups indexed by the simple roots. 

\gec  

\centerline{\bf 1.3. Moufang twin buildings}

\pec

The theory of {\it Moufang twin buildings } is the geometric side of the group
combinatorics above. A good account of their general theory is [A1]. It was initiated
in [T6], [T7] and {\it a posteriori } in [T1], whereas [RT2,3] deals with the special case
of twin trees. References for the classification problem are [MR] for
the unicity step, and then B. M\"uhlherr's work, in particular [M\"u] and his
forthcoming habilitationschrift. 

\pec

{\bf A. Twin buildings. } The definition of a {\it twin building } is
quite similar to that of a building in terms of $W$-distance [T7, A1 \S 2]. 

\pec 

{\bf Definition.---~} \it 
A {\rm twin building of type $(W,S)$ } consists in two
buildings  $({\cal I}_+,d_+)$ and $({\cal I}_-,d_-)$ of 
type $(W,S)$ endowed with a {\rm ($W$-)codistance}. By definition, the
latter is an application 
$d^*: ({\cal I}_+\times{\cal I}_-) \cup ({\cal I}_-\times{\cal I}_+) \rightarrow W$
satisfying the following conditions for each sign
$\epsilon$ and all chambers $x_\epsilon$ in ${\cal I}_\epsilon$ and 
$y_{-\epsilon}$, $y_{-\epsilon}'$ in ${\cal I}_{-\epsilon}$.

\pec

{\rm (TW1) } $d^*(y_{-\epsilon},x_\epsilon)=d^*(x_\epsilon,y_{-\epsilon})^{-1}$. 

{\rm (TW2) } If $d^*(x_\epsilon,y_{-\epsilon})=w$ and 
$d_{-\epsilon}(y_{-\epsilon},y_{-\epsilon}')=s \! \in \! S$ with $\ell(ws)<\ell(w)$, 
then $d^*(x_\epsilon,y_{-\epsilon}')=ws$.

{\rm (TW3) } If  $d^*(x_\epsilon,y_{-\epsilon})=w$ 
then for each  $s \! \in \! S$, 
there exists $z_{-\epsilon} \! \in \! {\cal I}_{-\epsilon}$ with 
$d_{-\epsilon}(y_{-\epsilon},z_{-\epsilon})=s$ and
$d^*(x_\epsilon,z_{-\epsilon})=ws$.
\rm\pec

From this definition can be derived two {\it opposition relations}. 
Two chambers are {\it opposite } if they are at codistance $1$. 
Given an apartment $\A_\epsilon$ of sign $\epsilon$, an {\it opposite }
of it is an apartment $A_{-\epsilon}$ such
that each chamber of $\A_\epsilon$ admits exactly one opposite in
$\A_{-\epsilon}$. In this situation, the same assertion is true after inversion of
signs. An apartment admits at most one opposite, and the set of 
apartments having an opposite forms an apartment system in the building [A1 \S 2]. 

\pec 

{\it Remark.---~} We defined two opposition relations, but we have to be careful with
them. Whereas an element of the apartment system defined above admits by definition 
exactly one opposite, a chamber admits many opposites in the building of opposite sign. 
In the sequel, the bold letter $\A$ will refer to a pair $(\A_+,\A_-)$ of opposite
apartments. 

\pec

The connection with group combinatorics is folklore [R\'e2, 2.6.4].  

\pec 

{\bf Proposition.---~} \it 
Let $G$ be a {\rm (TRD)}-group. Then the buildings associated to the two 
$BN$-pairs of $G$ are Moufang and belong to a twin building structure. 
In particular, two facets of opposite signs are always contained in a pair of opposite
apartments. At last, $G$ is transitive on pairs of opposite chambers. 
\cqfd\rm\pec 

So far, we explained how to derive two buildings from (TRD)-groups. 
The main consequences of  Moufang property will be described in the
next subsection. So what is left to do is to define the codistance. 
Whereas $W$-distance is deduced from Bruhat decomposition,
the $W$-codistance can be made completely
explicit thanks to the Birkhoff decomposition 1.2.B. Two chambers of
opposite signs $gB_+$ and $hB_-$ are at {\it codistance $w$} if and only if 
$g^{-1}h$ is in the Birkhoff class $B_+wB_-$. This definition does not
depend on the choice of $g$ and $h$ in their class. 

\pec

{\it Examples.---~} 1. The class of twin trees has been studied in full generality in
[RT2,3] where many properties are shown. For instance, the trees have to be at
least semi-homogeneous. The proofs are made simpler than in the general case thanks to the
use of an integral codistance which faithfully reflects the properties of the
$W$-codistance. 

2. For the groups ${\rm SL}_n(\K[t,t^{-1}])$ of 1.2.C,
the opposition relation can be formulated more concretely. Two
chambers $gB_+$ and $hB_-$ are opposite if and only if their
stabilizers -- Borel subgroups of opposite signs -- intersect along a
conjugate of $T$. Over $\fq$, this is equivalent to intersecting along a group of (minimal)
order $(q-1)^{n-1}$.  

\pec

{\bf B. Geometric description of the group action. } 
Suppose we are given a (TDR)-group $G$ and denote by $\A$ the
standard pair of opposite apartments in its twin building, that is 
$\{ wB_+ \}_{w \in W} \sqcup \{ wB_- \}_{w \in W}$. 
We want to give a geometric description of the $G$-action on its twin
building. 

\pec

{\it Kernel of the action.---~} One has (i) \it $H={\rm Fix}_G(\A_+)={\rm Fix}_G(\A_-)$.
\rm\pec

This result is the geometric formulation of [KP2, Corollary 3.4], it 
shows that the kernel of the action of $G$ on its twin building is contained in $H < G$. 

\pec

{\it Moufang property.---~} Roughly speaking, requiring the Moufang property to a
building is a way to make sure the latter admits a sufficiently large automorphism group,
with well-understood local actions. We just state the main consequence of it, which will
enable us to compute the number of chambers whose closure contains a given  panel $F$. This
integer will be referred to as the {\it thickness } at $F$ of the building. Let us start
with $F_s$ the positive panel of type $s$ in the closure of the standard chamber $c$. 
Then [Ron, (MO1) p.74]: \pec

(ii) \it The root group $U_{\alpha_s}$ fixes $\alpha_s$ and is normalized by $H$. It is
simply transitive on the chambers containing $F_s$ and $\neq c$. 
\rm\pec 

By homogeneity, point (ii) is true for every configuration panel-chamber-wall-root as
above, and {\it mutatis mutandis} everything remains true for the negative
building of $G$. 

\pec 

{\it Figure.---~} 

\centerline{\epsfysize=25mm $$\epsfbox{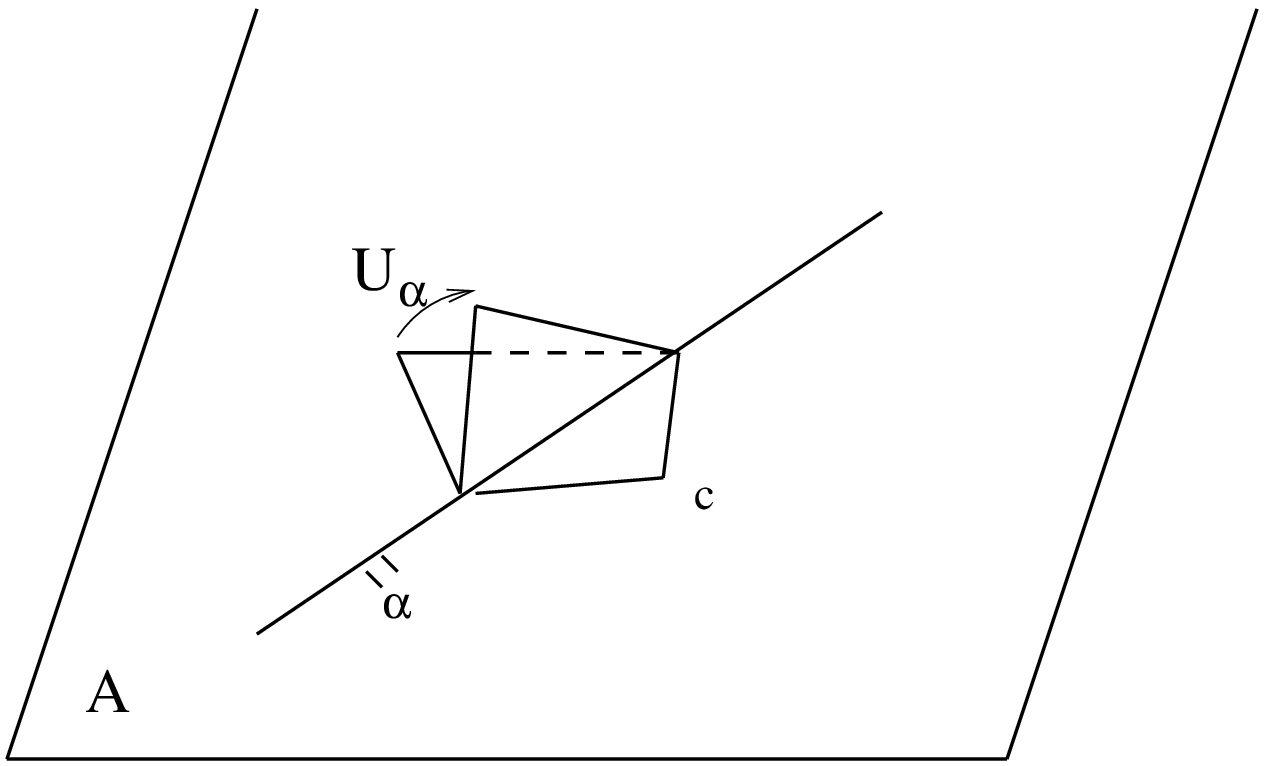}$$}

{\it The Weyl group as a subquotient.---~} The following assertion is just a
geometric formulation of axiom (TRD2). Consider the epimorphism $\nu: N
\twoheadrightarrow W$ with ${\rm Ker}(\nu)=H$. 

\pec

(iii) \it The group $N$ stabilizes $\A$; it permutes the $U_\alpha$'s by
$nU_\alpha n^{-1}=U_{\nu(n)\alpha}$, and the roots of $\A$ accordingly. 
\rm\gec

\centerline{\bf 1.4. A deeper use of geometry}

\pec 

In this section, we exploit different geometric notions to study in further detail 
(TRD)-groups. 

\pec 

{\bf A. Convexity and negative curvature. } To each of these notions corresponds a specific
geometric realization, the {\it conical } and the {\it metric } realization respectively. 
These geometries can be defined for a single building, they  basically differ by the way an
apartment is represented inside.  The representations of the whole building are obtained via
standard glueing techniques [D2, \S 10].

\pec

{\it The conical realization of a building.---~} An apartment of fixed sign is
represented here by the Tits cone of the Weyl group (1.1.B). 
For twin buildings, a pair of opposite apartments is represented by
the union of the Tits cone and its opposite in the ambient real vector
space. We can then use the obvious geometric opposition because it is the restriction
of the abstract opposition relation of facets. The point is that the shape of
simplicial cone for each facet fits in convexity arguments. 

\pec

{\it The metric realization of a building.---~} This realization was
defined by M. Davis and G. Moussong. Its main interest is that it is 
nonpositively curved. Actually, it satisfies the CAT(0) property: geodesic triangles are at least as thin as in the Euclidean plane. This
enables to apply the following [BrT1, 3.2] 

\pec 

{\bf Theorem (Bruhat-Tits fixed point theorem).---~} \it 
Every group of isometries of a {\rm CAT(0)}-space with a bounded orbit has a fixed point. 
\cqfd\rm\pec

This result generalizes a theorem applied by \'E. Cartan to Riemannian symmetric spaces to
prove conjugacy of maximal compact subgroups in Lie groups. The metric realization only represents
facets of spherical type. Indeed, consider the partially ordered set of spherical types $J
\subset S$. A chamber is represented by the cone over the barycentric subdivision of this
poset. G. Moussong [Mou] proved that this leads to a piecewise Euclidean cell complex
which is locally nonpositively curved. A simple connectivity criterion by M. Davis
[D1] proves the global CAT(0) property for Coxeter groups, and the use of retractions
proves it at the level of buildings [D2]. 

\pec

{\bf B. Balanced subsets. } In a single building, to fix a facet is the condition that
defines the so-called family of parabolic subgroups. In the twin situation, we can define
another family of subsets taking into account both signs. By taking fixators, it will also
give rise to an interesting family of subgroups. 

\pec

{\bf Definition.---~} \it Call {\rm balanced} a subset of a twin
building contained in a pair of opposite apartments, intersecting the
building of each sign and covered by a finite number of spherical facets. 
\rm\pec

Suppose we are given a balanced subset $\Omega$ and a pair of opposite apartments $\A$
containing it. According to the Bruhat decompositions, for each sign $\epsilon$ the
apartment $\A_\epsilon = \{ wB_\epsilon \}_{w \in W}$ is isomorphic to the Coxeter complex
of the Weyl group $W$.  We can then define interesting sets of roots associated to the
inclusion $\Omega \subset \A$, namely: \pec 

\centerline{$\Phi^u(\Omega) = \{ \alpha \in \Phi 
\mid \overline\alpha \supset \Omega \quad \partial\alpha \not\supset \Omega \}$, 
\quad $\Phi^m(\Omega) = \{ \alpha \in \Phi \mid \partial\alpha \supset \Omega \}$,} 

\pec

\centerline{and \quad 
$\Phi(\Omega) = \{ \alpha \in \Phi \mid \overline\alpha \supset \Omega \} 
= \Phi^u(\Omega)\sqcup\Phi^m(\Omega)$.}  

\pec

So as to work in an apartment of fixed sign, we shall use the terminology of separation and
strong separation. Using the Tits cone realization, we denote by
$\Omega_+$ (resp. $\Omega_-$) the subset $\Omega \cap \A_+$ (resp. the subset of opposites
of points in $\Omega \cap \A_-$).  This enables to work only in $\A_+$. The root $\alpha$ is
said to {\it separate } $\Omega_+$ from $\Omega_-$ if $\overline\alpha \supset \Omega_+$
while $-\overline\alpha \supset \Omega_-$; $\alpha$ {\it strongly separates} $\Omega_+$
from $\Omega_-$ if it separates $\Omega_+$ from $\Omega_-$ and if not both of $\Omega_+$
and $\Omega_-$ are contained in the wall $\partial\alpha$. Then $\Phi^u(\Omega)$ is the set
of roots strongly separating $\Omega_+$ from $\Omega_-$, $\Phi(\Omega)$ that of roots
separating $\Omega_+$ from $\Omega_-$. We will often omit the reference to $\Omega$ when it
is obvious. Here is an important lemma which precisely makes use of a convexity argument 
[R\'e2, 5.4.5]. 

\pec 

{\bf Lemma.---~ \it 
The sets of roots $\Phi^u$, $\Phi^m$ and $\Phi$ are finite. 
\rm\pec

{\it Idea of proof.---~} The set $\Phi^m$ is obviously stable under opposition, it is
finite since there is only a finite number of walls passing through a given spherical facet. 
Choose a point in each facet meeting $\Omega_+$ (resp. $\Omega_-$) and denote the
barycenter so obtained by $x_+$ (resp. $x_-$). The point $x_\pm$ is contained in a
spherical facet $F_\pm$. Connect the facets by a minimal gallery $\Gamma$, and denote by
$d_\pm$ the chamber whose closure contains $F_\pm$ and at maximal distance from the
chamber of $\Gamma$ whose closure contains $F_\pm$. By convexity, a root in $\Phi^u$ has to
contain $d_+$ but not $d_-$, so $\#\Phi^u$ is bounded by the length of a gallery connecting
these chambers. 

\pec 

{\it Figure.---~} 

\centerline{\epsfysize=25mm $$\epsfbox{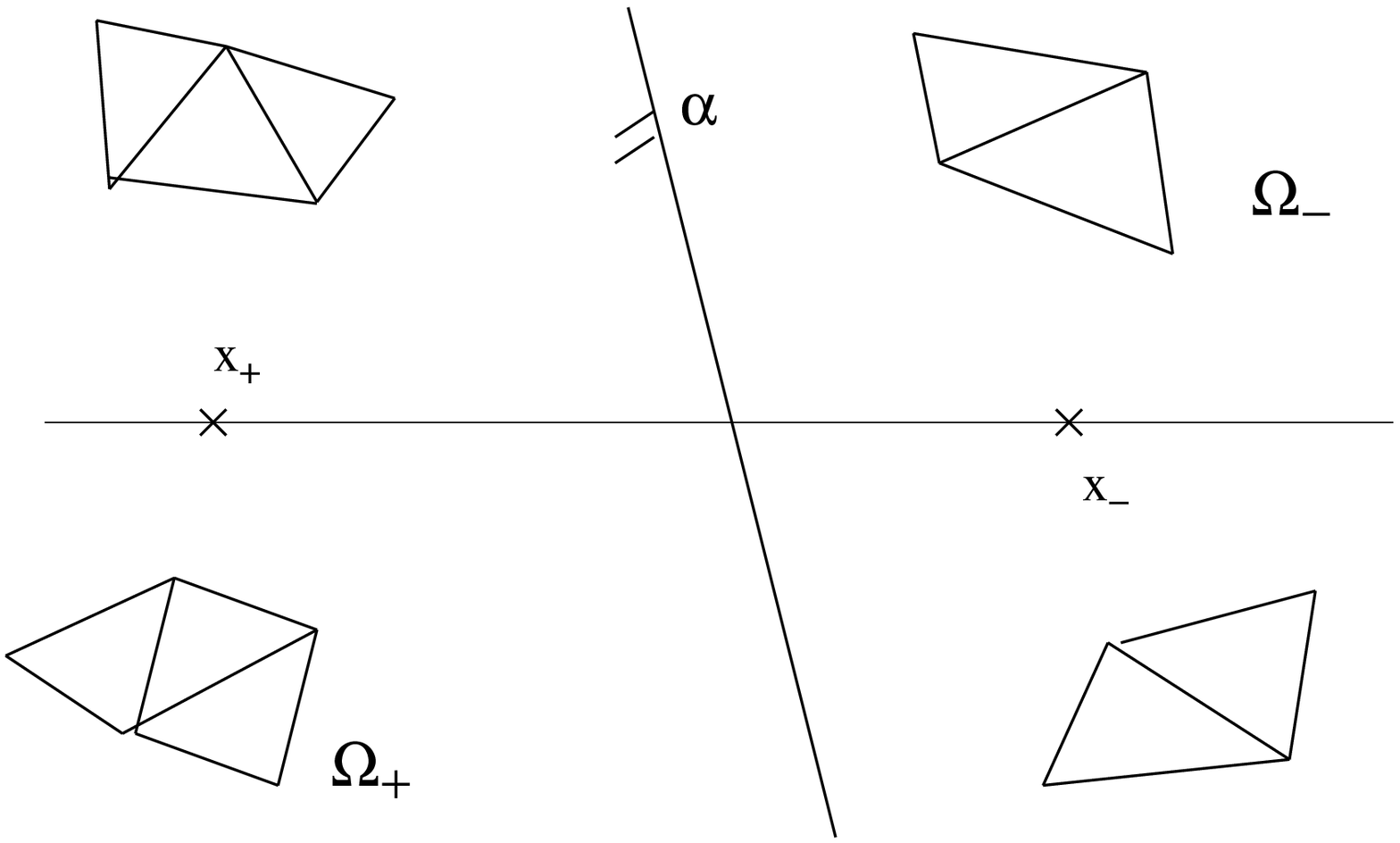}$$}

\cqfd\pec 

{\it Examples.---~} According to 1.3.A, a pair of spherical points of opposite
signs is a simple but fundamental example of balanced subset. Let us describe the sets
$\Phi^m$ and $\Phi^u$ in such a situation.  
Set $\Omega_\pm:= \{ x_\pm \} \subset \A_+$ so that $\Omega = \{ -x_-;x_+ \}$, and draw
the segment $[x_-;x_+]$ in the Tits cone. The set $\Phi^m$ is empty if $[x_-;x_+]$
intersects a chamber. The roots of $\Phi^u$ are the ones that strongly separate $x_+$
from $x_-$. 

\pec 

Assume first that each point lies in (the interior of) a chamber, which automatically
implies that $\Phi^m(\Omega)$ is empty. First, by transitivity of $G$ on the set of pairs of
opposite apartments (1.3.A), we may assume that $\A$ is the standard one. Then thanks to the
$W$-action, we can suppose that the positive point $x_+$ lies in the standard positive
chamber $c_+$; there is a $w$ in $W$ such that $x_-$ lies in $wc_+$. Choose a minimal
gallery between $c_+$ and $wc_+$. Then $\Phi^u$ is the set $\Phi_{w^{-1}}$ of the $\ell(w)$
positive roots whose wall is crossed by it. The first picture below is an example of type 
$\widetilde A_2$ with $\Phi^m=\varnothing$ and $\#\Phi^u=5$. 

\pec

Now we keep this $\widetilde A_2$ example, fix a wall $\partial\alpha$ and consider
pairs of points $\{ x_\pm \}$ with $[x_-;x_+] \subset \partial\alpha$. This implies that 
$\Phi^m$ will be always equal to $\{ \pm \alpha \}$. Still, as in the previous case, the
size of $\Phi^u$ will increase when so will the distance between $x_-$ and $x_+$. 
So the second picture below is an example of type $\widetilde A_2$ with $\#\Phi^m=2$ and
$\#\Phi^u=8$. 

\pec 

{\it Figure.---~} 

\centerline{\epsfysize=25mm $$\epsfbox{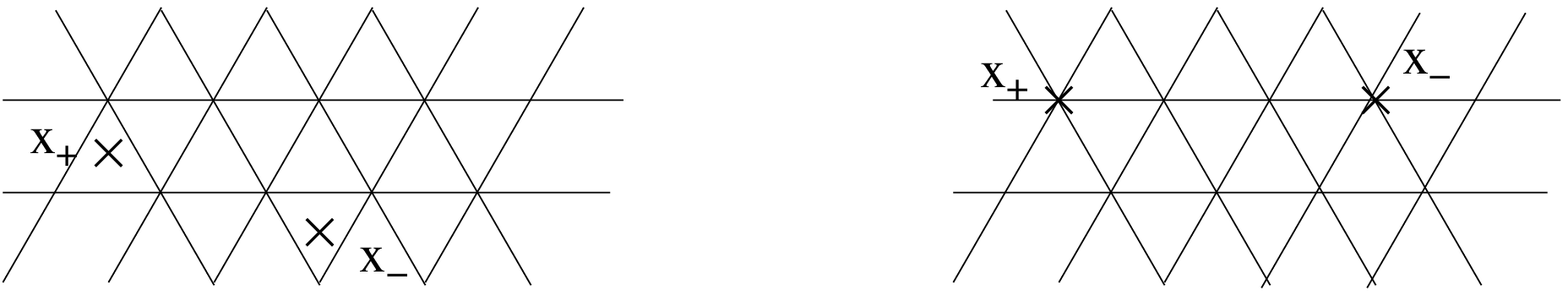}$$}

\centerline{\bf 1.5. Levi decompositions}

\pec

This section is dedicated to the statement of two kinds of
decomposition results, each being an abstract generalization of the
classical Levi decompositions. 

\pec

{\bf A. Levi decompositions for parabolic subgroups. } The first class of
subgroups distinguished in a group with a $BN$-pair is that of parabolic subgroups. The
existence of root groups in the (TRD) case suggests to look for more precise 
decompositions [R\'e2, 6.2]. 

\pec 

{\bf Theorem (Levi decomposition for parabolic subgroups).---~} \it
Let $F$ be a facet. Suppose $F$ is spherical or $G$ satisfies $\trdl$. 
Then for every choice of a pair $\A$ of opposite apartments containing $F$, the
corresponding parabolic subgroup ${\rm Fix}_G(F)$ admits a semi-direct product
decomposition 

\pec 

\centerline{${\rm Fix}_G(F) = M(F,\A) \ltimes U(F,\A)$.}

\pec 

The group $M(F,\A)$ is the fixator of $F \cup -F$, where $-F$ is the opposite of $F$ in
$\A$; it is generated by $H$ and the root groups $U_\alpha$ with $\partial\alpha \supset
F$. It satisfies the {\rm (TRD) } axioms, and the $\trdl$ refinement if so does $G$. 
The group $U(F,\A)$ only depends on $F$, it is the normal closure in
${\rm Fix}_G(F)$ of the root groups $U_\alpha$ with $\alpha \supset F$. 
\cqfd\rm\pec

{\it Remarks.---~} 1. It is understood in the statement that all the roots considered
above are subsets of $\A$. By transitivity of $G$ on pairs of opposite apartments, the root
groups are conjugates of that in the (TRD) axioms. 

2. In the case of an infinite Weyl group, Levi decompositions show that the action of a
(TRD)-group on a single building cannot be discrete since the fixator of a point contains
infinitely many root groups. 

\pec 

{\bf B. Levi decompositions for small subgroups. } As mentioned in 1.4.B, the twin
situation leads to the definition of another class of interesting subgroups. 

\pec

{\bf Definition.---~} \it A subgroup of $G$ is called {\rm small }
if it fixes a balanced subset. 
\rm\pec

To describe the combinatorial structure of small subgroups, we need to require the 
additional (NILP) axiom. It is not really useful to state it here
explicitly, because our basic object of study is the class of Kac-Moody groups which
satisfy it. We have  [R\'e2, 6.4]: \pec

{\bf Theorem (Levi decomposition for small subgroups).---~} \it 
Suppose $G$ is a $\trdl$-group satisfying {\rm (NILP)}, and let $\Omega$ be a balanced
subset.  Then for every choice of pair of opposite apartments $\A$ containing $\Omega$, the
small subgroup ${\rm Fix}_G(\Omega)$ admits a semi-direct product decomposition 

\pec 

\centerline{${\rm Fix}_G(\Omega) = M(\Omega,\A) \ltimes U(\Omega,\A)$,}

\pec 

where both factors are closely related to the geometry of $\Omega$
(and of $\A$, w.r.t. which all roots are defined). 
Namely, $M(\Omega,\A)$ satisfies the $\trdl$ axioms for the root
groups $U_\alpha$ with $\partial\alpha \supset \Omega$, 
that is $\alpha \! \in \! \Phi^m(\Omega)$; 
and $U(\Omega,\A)$ is in bijection with the set-theoretic product of the root
groups $U_\alpha$ with $\alpha \! \in \! \Phi^u(\Omega)$ for any given order. 
\cqfd\rm\pec

{\it Remark.---~} In the Kac-Moody case, it can be shown that $U(\Omega,\A)$ only
depends on $\Omega$ and not on $\A$, so that we will use the notation $U(\Omega)$ instead. 

\pec 

{\it Examples.---~} 1. As a very special case, consider in the standard pair of
opposite apartments the balanced subset $\Omega:= \{ c_+;wc_- \}$, where $c_\pm$ are the
standard chambers. Then, thanks to the first case of example 1.4.B, we
see that ${\rm Fix}_G(\Omega) = T \ltimes U_w$, with $U_w$ as defined in 1.2.B.
Set-theoretically, the latter group is in particular a product of $\ell(w)$ root groups. 

2. Another special case is provided by the datum of two opposite spherical facets $\pm F$ in
$\A$. Then $U(F \cup -F) = \{ 1 \}$ and the fixator of the pair is at the same time the
$M$-factor of the Levi decomposition above and in the previous sense for the parabolic
subgroups ${\rm Fix}_G(F)$ and ${\rm Fix}_G(-F)$. 

\pec

{\bf C. } {\it Further examples.---~} The examples below appear in the two simplest cases
of special linear groups over Laurent polynomials (1.2.C). They will illustrate both
kinds of decomposition. Examples of hyperbolic twin buildings will be treated in \S 4.1. 

\pec 

{\it The ${\rm SL}_2(\fq[t,t^{-1}])$ case.---~} The Bruhat-Tits buildings involved are
both isomorphic to the homogeneous tree $T_{q+1}$ of valency $q+1$. A pair of opposite
apartments is represented by two parallel real lines divided by the integers. 

\pec

Let us consider first the fixator of a single point. If this point is
in an open edge, then its fixator is a Borel subgroup isomorphic to
the group ${\rm SL}_2\pmatrix{\fq[t] & \fq[t] \cr t\fq[t] & \fq[t]}$. If the point is a
vertex, then it is isomorphic to the Nagao lattice ${\rm SL}_2(\fq[t])$. 

\pec

Let us consider now pairs of points $\{ x_+;-x_- \}$ of opposite signs. 
We use the notations of the examples of 1.4.B. If the points are opposite vertices, then
$\Phi^u$ is empty and $\Phi^m$ consists of two opposite roots: the corresponding small
subgroup is equal to the $M$-factor, a finite group of Lie type and rank one over $\fq$. If
the points are non opposite middles of edges $e_+$ and $-e_-$, and if the edges $e_+$
and $e_-$ in the positive straight line are at $W$-distance $w$, then the corresponding
small fixator is the semi-direct product of $T$ (isomorphic to $\fq^\times$) by a
commutative group isomorphic to the additive group of the $\fq$-vector space $\fq^{\ell(w)}$. 

\pec 

{\it Figure.---~} 

\centerline{\epsfysize=30mm $$\epsfbox{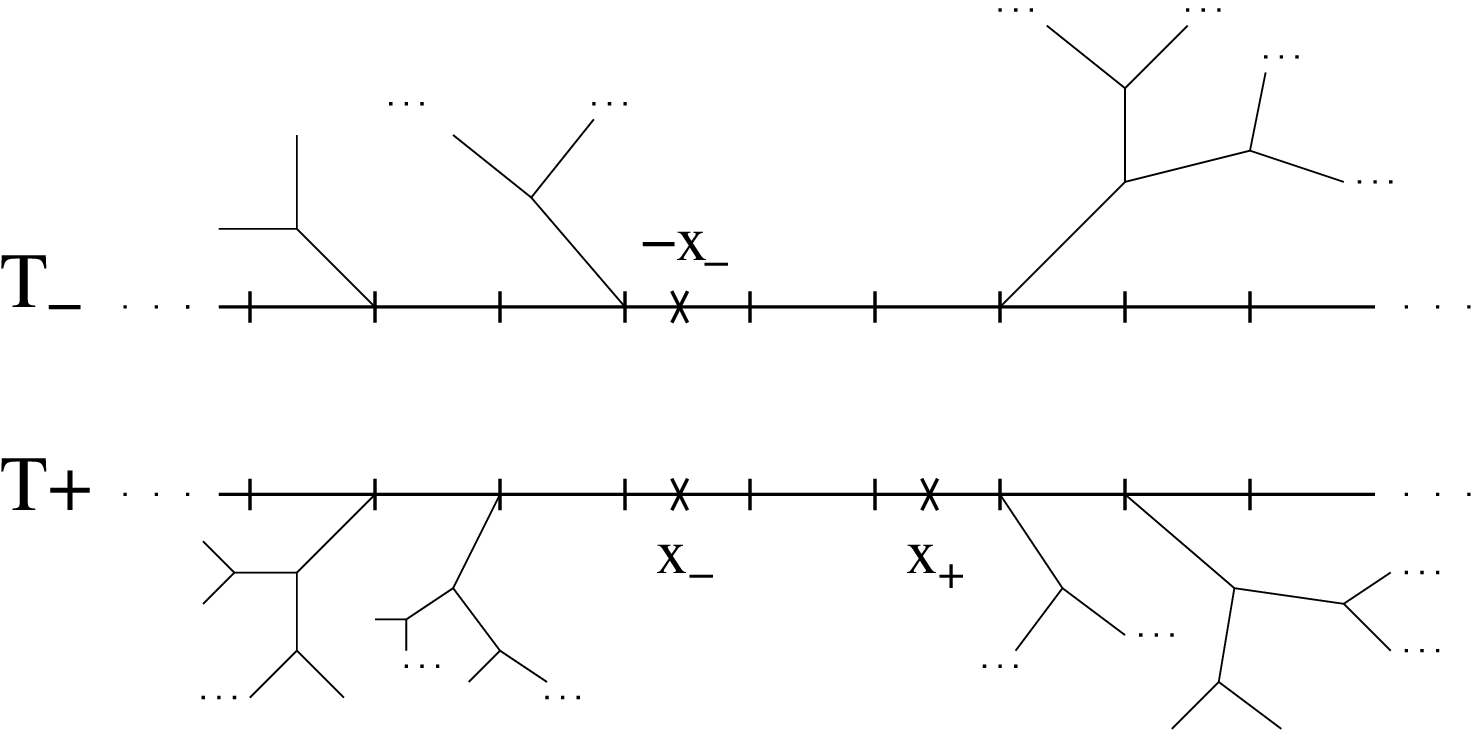}$$}

{\it The ${\rm SL}_3(\fq[t,t^{-1}])$ case.---~} The case of single points is
quite similar. For instance, let us fix a vertex $x$. The set $L$ of chambers whose
closure contains $x$ can naturally be seen as the building of ${\rm SL}_3(\fq)$. The
fixator $P$ of $x$ is isomorphic to ${\rm SL}_3(\fq[t])$. The $M$-factor in the Levi
decomposition of $P$ is a finite group of Lie type $A_2$ over $\fq$ and its action on the
small building $L$ is the natural one. The infinite $U$-factor fixes $L$ pointwise. 

\pec

Suppose now we are in the second case of example 1.4.B: the pair of points $\{ x_+;-x_- \}$
lies in a wall $\partial\alpha$. We already now that $\Phi^m$ contains 
the two opposite roots $\pm \alpha$, so that the $M$-factor of Fix($\{ x_+;-x_-\}$) is a
rank one finite group of Lie type over $\fq$. If we assume $\alpha$ to be the first simple
root of ${\rm SL}_3(\fq)$, then $M$ is the group of matrices 
$\pmatrix{X & 0 \cr 0 & {\rm det}X^{-1}}$, with $X \! \in \! {\rm GL}_2(\fq)$. 
The $U$-factor is a (commutative) $p$-group, with $p = {\rm char }(\fq)$. It is isomorphic
to the group of matrices $\pmatrix{1 & 0 & P \cr 0 & 1 & Q \cr 0 & 0 & 1}$, where $P$ and $Q$
are polynomials of $\fq[t]$ of degree bounded the number of vertices between $x_-$ and
$x_+$. 

\vfill\eject

\centerline{\bf 2. SPLIT KAC-MOODY GROUPS}

\gec

In 1987, J. Tits defined group functors whose values over fields we will call 
{\it (split) Kac-Moody groups}. In \S 2.1, we give the definition of
these groups, which involves some basic facts from Kac-Moody algebras
to be recalled. Section 2.2 relates Kac-Moody groups to the abstract
combinatorics previously studied, providing afterwards some more
information. The three last sections deal then with more specific
properties. In \S 2.3 is presented the adjoint representation and \S
2.4 combine it with Levi decompositions to endow each small subgroup with
a structure of algebraic group. At last, \S 2.5 presents an argument
occuring several times afterwards. It is used for example to prove a conjugacy theorem for
Cartan subgroups. 

\gec

\centerline{\bf 2.1. Tits functors and Kac-Moody groups}

\pec

We sum up here the step-by-step definition of split Kac-Moody groups by J. Tits. It is also
the opportunity to indicate how the notions defined in \S 1.1 for arbitrary Coxeter groups
arise in the Kac-Moody context. All the material here is contained in [T4, \S 2]. 

\pec

{\bf A. Definition data. Functorial requirements. } The data needed to define a Tits group
functor are quite similar to that needed to define a Chevalley-Demazure group scheme. 

\pec 

{\bf Definition.---~} \it 
{\rm (i) } A {\rm generalized Cartan matrix } is an integral matrix
$A=[A_{s,t}]_{s,t \in S}$ satisfying: $A_{s,s}=2$, $A_{s,t} \leq 0$
when $s \neq t$ and $A_{s,t}=0 \Leftrightarrow A_{t,s}=0$. 

{\rm (ii) } A {\rm Kac-Moody root datum } is a $5$-tuple 
$\bigl(S, A, \Lambda, (c_s)_{s \in S}, (h_s)_{s \in S} \bigr)$, where $A$ is a 
generalized Cartan matrix indexed by the finite set $S$, $\Lambda$ is a free $\Z$-module 
(whose $\Z$-dual is denoted by $\tch{\Lambda}$). Besides, the elements $c_s$ of $\Lambda$
and $h_s$ of $\tch{\Lambda}$ are required to satisfy $\langle c_s \mid h_t \rangle
=A_{ts}$ for all $s$ and $t$ in $S$. 

\rm\pec 

Thanks to these elementary objects, it is possible to settle
abstract functorial requirements generalizing the properties of functors of points of
Chevalley-Demazure group schemes [T4, p.545]. The point is that {\it over fields} a
problem so defined does admit a unique solution, concretely provided by a group
presentation [T4, Theorem 1]. This is the viewpoint we will adopt in the next subsection and
for the rest of our study. Let us first define further objects only depending on a Cartan
generalized matrix. 

\pec 

{\bf Definition.---~} \it 
{\rm (i) } The {\rm Weyl group } of a generalized Cartan matrix $A$ is the Coxeter group  

\pec 

\centerline{$W = \langle s \!\in \! S \mid (st)^{M_{st}}=1$ whenever $M_{st} < \infty 
\rangle$.}

\pec

with  $M_{st} = 2, 3, 4, 6$ or $\infty$ according to whether
$A_{st}A_{ts} = 0, 1, 2, 3$ or is $\geq 3$. 

{\rm (ii) } The {\rm root lattice } of $A$ is the free $\Z$-module 
$\displaystyle Q:= \bigoplus_{s \in S} \Z a_s$ over the symbols
$a_s$, $s \! \in \! S$. 
We also set $\displaystyle Q_+:=\bigoplus_{s \in S} \N a_s$. 
\rm\pec 

The group $W$ operates on $Q$ via the formulas $s.a_t=a_t - A_{st}a_s$. 
This allows to define the {\it root system } $\Delta^{re}$ of $A$ by 
$\Delta^{re}:= W.\{ a_s \}_{s \in S}$. This set has all the abstract
properties of the set of real roots of a Kac-Moody algebra. In
particular, one has $\Delta^{re} = \Delta^{re}_+ \sqcup
\Delta^{re}_-$, where $\Delta^{re}_\epsilon:= \epsilon Q_+ \cap \Delta^{re}$. 
Let us explain now how to recover the notions of \S 1.1. Recall
that given a Coxeter matrix $M$, we defined its cosine matrix so as to introduce the Tits
cone on which acts the group. The starting point here is the generalized Cartan matrix
which plays the role of the cosine matrix. This allows to define the Tits cone as in 1.1.B,
and all the objects related to it. Instead of $V$, the real vector space we use here is 
$Q_{\pr}:= Q \otimes_{\pz} \R$. This way, we get a bijection between $\Delta^{re}$ and the
root system of the Weyl group as defined in 1.1.A, for which the notions of prenilpotence
and intervals coincide. General references for abstract infinite root systems are
[H2] and [Ba]. Note that in the Kac-Moody setting, linear combinations of roots
giving roots only involve integral coefficients. This leads to the following

\pec 

{\bf Convention.---~} \it 
In Kac-Moody theory, the root system $\Delta^{re}$ will play the role of the root system of
the Weyl group of $A$, and the roots will be denoted by latin letters. 
\rm\pec

{\bf B. Lie algebras. Universal envelopping algebras. Z-forms. } 
As already said, a more constructive approach consists in defining groups by generators and
relations. The idea is to mimic a theorem by Steinberg giving a presentation of a split
simply connected semisimple algebraic group over a field. So far, so good. The point is that
the relations involving root groups are not easily available. One has to reproduce some
kind of Chevalley's construction in the Kac-Moody context. This imply to introduce a quite
long series of \og tangent\fg objects such as Lie algebras, universal envelopping algebras
and at last $\Z$-forms of them.  

\pec 

{\it Kac-Moody algebras.---~} 
From now on, we suppose we are given a Kac-Moody root datum 
${\cal D} = \bigl( S, A, \Lambda, (c_s)_{s \in S}, (h_s)_{s \in S} \bigr)$. 

\pec

{\bf Definition.---~} \it
The {\rm Kac-Moody algebra } associated to ${\cal D}$ is the $\C$-Lie
algebra $\g_{\cal D}$ generated by $\g_0:= \tch{\Lambda} \otimes_{\pz} \C$ and 
the sets $\{ e_s \}_{s \in S}$ and $\{ f_s \}_{s \in S}$, submitted to
the following relations.

\pec

\centerline{$[h, e_s]= \langle c_s, h \rangle e_s \quad  
{\it and }  \quad [h, f_s]= -\langle c_s, h \rangle f_s$ \ 
{\it for $h \! \in \! \g_0$}, \quad $[\g_0,\g_0]=0$,}

\pec 

\centerline{$[e_s,f_s]= -h_s \otimes 1$, \quad $[e_s,f_t]=0$ 
\ {\it for $s \neq t$ in $S$},}

\pec

\centerline{$(\ad e_s)^{-A_{st}+1}e_t = (\ad f_s)^{-A_{st}+1}f_t = 0$
\quad \hbox{\it -- \quad Serre relations.}}
\rm\pec

This Lie algebra admits an abstract $Q$-gradation by declaring the
element $e_s$ (resp. $f_s$) of degree $a_s$ (resp. $-a_s$) and the
elements of $\g_0$ of degree $0$. Each degree for which the
corresponding subspace is nontrivial is called a {\it root } of
$\g_{\cal D}$. 

\pec 

For each $s \! \in \! S$, the derivations $\ad e_s$ and $\ad f_s$ are
locally nilpotent, which enables to define the automorphisms
$\exp(\ad e_s)$ and $\exp(\ad f_s)$ of $\gd$. 
The automorphisms $\exp(\ad e_s).\exp(\ad f_s).\exp(\ad e_s)$ and 
$\exp(\ad f_s).\exp(\ad e_s).\exp(\ad f_s)$ are equal and will be denoted by $s^*$. We are
interested in $s^*$ because if we define $W^*$ to be the group of automorphisms of $\gd$
generated by the $s^*$'s, then 
$s^* \mapsto s$ lifts to a homomorphism $\nu: W^* \twoheadrightarrow W$. 
Considering the $W^*$-images of the spaces $\C e_s$ and $\C f_s$ makes  
$\Delta^{re}$ appear as a subset of roots of $\gd$. A basic fact about real roots $a$ in
Kac-Moody theory is that the corresponding homogeneous subspaces $\g_a \subset \gd$ are
one-dimensional, so that the $W$-action on $\Delta^{re}$ is lifted by the $W^*$-action on the
homogeneous spaces $\g_a$, $a \! \in \! \Delta^{re}$. 
Besides, for $s \! \in \! S$ and  $w^* \! \in \! W^*$, the pair of opposite elements 
$w^*(\{ \pm e_s \})$ only depends on the root $a:=\nu(w^*)a_s$; we denote it by 
$\{ \pm e_a \}$. 

\pec

{\it Divided powers and $\Z$-forms.---~} We want now to construct a
$\Z$-form of the universal envelopping algebra $\ugd$ in terms of divided powers. Subrings
of this ring will be used for algebraic differential calculus on some subgroups of Kac-Moody
groups. For each $u \! \in \! \ugd$ and each $n \! \in \! \N$, $u^{[n]}$ will denote the
divided power $(n!)^{-1}u^n$ and $\pmatrix{u \cr n \cr}$ will denote
$(n!)^{-1}u(u-1)...(u-n+1)$.  For each  $s \! \in \! S$, 
${\cal U}_{\{ s \} }$ (resp. ${\cal U}_{\{ -s \} }$)  is the subring 
$\sum_{n \in \pn} \Z e_s^{[n]}$ (resp. $\sum_{n \in \pn} \Z f_s^{[n]}$) of
$\ugd$. We denote by $\uo$ the subring of $\ugd$ generated by the (degree 0) elements of the
form $\pmatrix{h \cr n \cr}$, with $h \! \in \! \tch{\Lambda}$ and $n \! \in \! \N$. 

\pec 

{\bf Definition.---~} \it 
We denote by $\ud$ the subring of $\ugd$ generated by $\uo$ and the 
${\cal U}_{\{ s \} }$ and ${\cal U}_{\{ -s \} }$ for $s$ in $S$. 
It contains the ideal $\ud^+:= \gd . \ugd \cap \ud$. 
\rm\pec 

The point about $\ud$ is the following 

\pec 

{\bf Proposition.---~} \it 
The ring $\ud$ is a {\rm $\Z$-form } of $\ugd$ i.e., the natural
map $\ud \otimes_{\pz} \C \to \ugd$ is a bijection. 
\rm\pec

{\it Remark.---~} In the finite-dimensional case -- when $\gd$ is a semisimple Lie 
algebra, $\ud$ is the ring used to define Chevalley groups by means of formal exponentials
[Hu1, VII]. 

\pec 

{\bf C. Generators and Steinberg relations. Split Kac-Moody groups. }
We can now start to work with groups. 

\pec 

{\it First step: root groups and integrality result.---~} Let us consider a root $c$ in
$\Delta^{re}$.  We denote by ${\goth U}_c$ the $\Z$-scheme isomorphic to the additive
one-dimensional group scheme ${\bf G}_a$ and whose Lie algebra is the $\Z$-module generated
by $\{ \pm e_c \}$. For each root $c$, the choice of a sign defines an isomorphism 
${\bf G}_a \iso {\goth U}_c$. These choices are analogues of the ones defining an \og
\'epinglage\fg in the finite-dimensional case. 
We denote for short ${\goth U}_s$ (resp. ${\goth U}_{-s}$) 
the $\Z$-group scheme ${\goth U}_{a_s}$ (resp. ${\goth U}_{-a_s}$), and
$x_s$ (resp. $x_{-s}$) the isomorphism ${\bf G}_a \iso {\goth U}_s$ 
(resp. ${\bf G}_a \iso {\goth U}_{-s}$) induced by the choice of $e_s$ 
(resp. $f_s$) in $\{ \pm e_s \}$ (resp. $\{ \pm f_s \}$). 
Let us consider now $\{ a;b \}$ a prenilpotent pair of roots. The
direct sum of the homogeneous subspaces $\g_c$ of degree 
$c \! \in \! [ a;b ]_{\hbox{\sevenrm lin}}$ is a nilpotent Lie subalgebra 
$\g_{[ a;b ]_{\hbox{\sevenrm lin}}}$ of $\gd$, which defines a unique unipotent complex
algebraic group $U_{[ a;b ]_{\hbox{\sevenrm lin}}}$, by means of formal exponentials. 

\pec 

{\bf Proposition.---~} \it 
There exists a unique $\Z$-group scheme ${\goth U}_{[ a;b ]_{\hbox{\sevenrm lin}}}$
containing the $\Z$-schemes ${\goth U}_c$ for  $c \! \in \! 
[ a;b ]_{\hbox{\sevenrm lin}}$, whose value over $\C$ is the complex algebraic group 
$U_{[ a;b ]_{\hbox{\sevenrm lin}}}$
with Lie algebra $\g_{[ a;b ]_{\hbox{\sevenrm lin}}}$. Moreover, the product map 

\pec

\centerline{$\displaystyle \prod_{c \in [ a;b ]_{\hbox{\sevenrm lin}}} 
{\goth U}_c \rightarrow {\goth U}_{[ a;b ]_{\hbox{\sevenrm lin}}}$}

\pec 

is an isomorphism of $\Z$-schemes for every order on $[ a;b ]_{\hbox{\sevenrm lin}}$. 
\rm\pec

{\it Remark.---~} This result is an abstract generalization of an 
integrality result of Chevalley [Sp, Proposition 9.2.5] concerning the commutation
constants between root groups of split semisimple algebraic groups. The
integrality here is expressed by the existence of a $\Z$-structure
extending the $\C$-structure of the algebraic group. We will make these $\Z$-group schemes
a bit more explicit when defining the adjoint representation (in 2.3.A). 

\pec 

{\it Second step: Steinberg functors.---~} Now we see the $\Z$-schemes
${\goth U}_{[ a;b ]_{\hbox{\sevenrm lin}}}$ just as group functors over rings. We will
amalgamate them so as to obtain a big group functor containing all the relations between the
root groups. This is the step of the procedure where algebraic structures are lost. 
The relations $c \! \in \! [ a;b ]_{\hbox{\sevenrm lin}}$ give rise to an inductive system
of group functors with injective transition maps 
${\goth U}_c(R)  \hookrightarrow {\goth U}_{[ a;b ]_{\hbox{\sevenrm lin}}}(R)$ 
over every ring $R$. 

\pec 

{\bf Definition.---~} \it 
The {\rm Steinberg functor } -- denoted by $\st$ -- is the limit of the inductive system
above. 
\rm\pec 

This group functor only depends on the generalized Cartan matrix $A$, not on the whole datum
${\cal D}$. Besides, the $W^*$-action on $\gd$ extends to a $W^*$-action on $\st$: for each
$w$ in $W$, we will still denote by $w^*$ the corresponding automorphism. 

\pec 

{\it Third step: the tori and the other relations.---~} What is left to do at this level
is to introduce a torus, a split one by analogy with Chevalley groups. In fact, it is
determined by the lattice $\Lambda$ of the Kac-Moody root datum ${\cal D}$. As in the
previous step, we are only interested in group functors. Denote by ${\cal T}_\Lambda$ the
functor of points 
$\Hom_{\pz -{\hbox{\sevenrm alg}}}(\Z[\Lambda],-)
=\Hom_{{\hbox{\sevenrm groups}}}(\Lambda,-^\times)$ 
defined over the \og category\fg $\Z$-alg of rings. 
For a ring $R$, $r$ an invertible element in it and $h$ an element of
$\tch{\Lambda}$, the notation $r^h$ corresponds to the element 
$\lambda \mapsto r^{\langle \lambda \mid h \rangle}$ of ${\cal T}_\Lambda(R)$. 
We will also denote $r^{\langle \lambda \mid h \rangle}$ by $\lambda(r^h)$. 
This way we can see elements of $\Lambda$ as characters of ${\cal T}_\Lambda$.  We finally
set $\tilde s(r):=x_s(r)x_{-s}(r^{-1})x_s(r)$ for each invertible element $r$ in a ring $R$,
and $\tilde s:= \tilde s(1)$. 

\pec

{\bf Definition.---~} \it
For each ring $R$, we define the group $\gkm(R)$ as the quotient of
the free product $\st (R) * {\cal T}_{\Lambda}(R)$ by the following
relations. 

\pec 

\centerline{$tx_s(q)t^{-1}=x_s(c_s(t)q); \quad \tilde s(r) t \tilde s(r)^{-1}=s(t)$;}

\pec

\centerline{$\tilde s(r^{-1})=\tilde s r^{h_s}; \quad 
\tilde s u {\tilde s}^{-1} = s^*(u)$;}

\pec 

with $s \! \in \!S,  c \! \in \! \Delta^{re},  t \! \in \! {\cal T}_{\Lambda}(R),  
q \! \in \! R, r \! \in \! R^\times, u \! \in \! {\goth U}_c (R)$. 

\pec

We define this way a group functor denoted by $\gkm$ and called the 
{\rm Tits functor } associated to ${\cal D}$. 
The value of this functor on the field $\K$ will be called the {\rm split
Kac-Moody group } associated to ${\cal D}$ and defined over $\K$. 
\rm\pec 

{\it Remark.---~} This is the place where we use the  duality relations
between $\Lambda$ and $\tch{\Lambda}$ required in the definition of ${\cal D}$. 
This duality occurs in the expression $c_s(t)$, and allows to define a $W$-action on 
$\Lambda$ by sending $s$ to the involution of $\Lambda$: $\lambda \mapsto \lambda - \langle \lambda \mid h_s \rangle c_s$. This yields a $W$-action
on ${\cal T}_{\Lambda}$ thanks to which the expression $s(t)$ makes sense in the defining
relations. The expression $s^*(u)$ comes from the $W^*$-action on the Steinberg functor. 

\pec 

{\it Examples.---~} 1. As a first example, we can consider a Kac-Moody root datum
where $A$ is a Cartan matrix and $\tch{\Lambda}$ is the $\Z$-module freely generated by
the $h_s$'s. Then what we get over fields is the functor of points of the simply connected Chevalley
scheme corresponding to $A$. This can be derived from the theorem by Steinberg alluded to
in 2.1.A [Sp, Theorem 9.4.3]. One cannot expect better since as already explained,
algebraic structures are lost during the amalgam step defining $\st$. 

2. More generally, the case where $\Lambda$ (resp. $\tch{\Lambda}$) is freely generated
by the $c_s$'s (resp. $h_s$'s) will be referred to as the {\it adjoint } (resp. {\it simply
connected}) case. 

3. The groups ${\rm SL}_n(\K[t,t^{-1}])$ of 1.2.C can be seen as split Kac-Moody groups over
$\K$. The generalized Cartan matrix is the matrix $\widetilde A_{n-1}$ indexed by $\{ 0;
1; ... n-1 \}$, characterized for $n \geq 3$ by $A_{ii}=2$, $A_{ij}=-1$ for $i$ and $j$
consecutive modulo $n$ and $A_{ij}=0$ elsewhere. In rank 2, $\widetilde A_2$ is just 
$\pmatrix{2 & -2 \cr -2 & 2}$. The $\Z$-module $\Lambda$ (resp. $\tch{\Lambda}$) is then
the free $\Z$-module of rank $n-1$ generated by the simple roots $c_1$, ... $c_{n-1}$ of the
standard Cartan subgroup $T$ of ${\rm SL}_n(\K)$ (resp. the corresponding one-parameter
multiplicative subgroups). To be complete, one has to add the character 
$c_0:=-\sum_{1 \leq i \leq n-1} c_i$ 
(resp. the cocharacter $h_0:=-\sum_{1 \leq i \leq n-1} h_i$). 
These examples show in
particular that neither $\Lambda$ nor $\tch{\Lambda}$ need to contain the $c_s$'s or the
$h_s$'s as free families. 

\gec 

\centerline{\bf 2.2. Combinatorics of a Kac-Moody group}

\pec 

This section states the connection with the abstract group
combinatorics presented in \S 1, then stresses some specific
properties (often due to splitness). We are given here a split Kac-Moody group $G$ 
defined over $\K$ which acts on its twin building $(\ip,\im,d^*)$.

\pec

{\bf A. Main statement. } As the (TRD) axioms were adjusted to the Kac-Moody situation, it
is not surprising to get the following result [R\'e2, 8.4.1]. 

\pec 

{\bf Proposition.---~} \it 
Let $G$ be a Kac-Moody group, namely the value of a Tits group functor
over a field $\K$. Then $G$ is a $\trdl$-group for the root groups of 
its definition, which also satisfies also the axiom {\rm (NILP) } for Levi decomposition
of small subgroups {\rm (1.5.B)}. The standard Cartan subgroup $T$ plays the role of the
normalizer of the $U_\alpha$'s. 
\rm\pec 

{\it Sketch of proof.---~} In the defining relations of a Kac-Moody
group (2.1.C), everything seems to be done so as to conform to the (TRD)
axioms. Still, we have to be a bit more careful with the non
degeneracy requirements. Indeed, for the first half of axiom (TRD0)
and axiom (TRD3), we have to use the adjoint representation 
defined later. This point is harmless since the definition of this
representation does not need group combinatorics. 
\cqfd\pec 

{\it Remark.---~} This result can be used to prove the unicity
over fields of Kac-Moody groups. Actually, the original group
combinatorics involved in the proof by J. Tits [T4, \S 5] was more general. 

\pec 

{\bf B. Specific properties. } As Tits functors generalize Chevalley-Demazure schemes, one
can expect further properties for at least two reasons. The first one
is the \og algebraicity\fg of Kac-Moody groups i.e., the fact that the precise defining
relations are analogues of that for algebraic groups. The second reason is splitness. 

\pec 

{\it \og Algebraicity\fg of Tits functors.---~} First, we know that the biggest group
normalizing all the root groups is the standard Cartan subgroup $T$: this is a purely
combinatorial consequence of the (TRD) axioms due to the fact that in $BN$-pairs, Borel
subgroups coincide with their normalizer. Besides we have a completely explicit description
of the $T$-action on root groups -- by characters. This leads to the following result 
[R\'e2, 10.1.3]. 

\pec 

{\bf Proposition.---~} \it 
If the groundfield $\K$ has more than $3$ elements, then the fixed points under $T$ in the
buildings are exactly the points in $\A$ $($the standard pair of opposite apartments$)$. 
\cqfd\rm\pec

{\it Remark.---~} What this proposition says is that this set is
not bigger than $\A$. This result has to be related to another one,
where the same hypothesis $\# \K > 3$ is needed: in this case, $N$ is
exactly the normalizer of $T$ in $G$ [R\'e2, 8.4.1]. 

\pec 

Another algebraic-like consequence is about Levi decompositions and residues in buildings.
Consider the standard facet $\pm F_J$ of spherical type $J \subset S$ and sign $\pm$. Then
the $M$-factor (see 1.5.B, example 2) is the Kac-Moody group corresponding to the datum
where $A$ is restricted to $J \times J$, and only the $c_s$'s and $h_s$'s with $s \!
\in \! J$ are kept. This is a root datum abstractly defining points of a reductive group.
Moreover, the union of the sets of chambers whose closure contains $\pm F_J$ -- the {\it
residue } of $\pm F_J$ -- is the twin building of the reductive
group, on which it operates naturally. The
generalization to arbitrary facets is straightforward. 

\pec

{\it Remark.---~} The analogy with reductive groups can fail on some points. Consider
for instance the case of the centralizer of the standard Cartan
subgroup $T$. It is classical in the algebraic setting that the
centralizer of a Cartan subgroup is not bigger than the subgroup itself. 
Consider once again the groups ${\rm SL}_n(\K[t,t^{-1}])$ of 1.2.C. 
Then the centralizer of $T$ -- the standard Cartan subgroup of 
${\rm SL}_n(\K)$ -- is infinite-dimensional since it is the group of
determinant 1 diagonal matrices with monomial coefficients. 

\pec 

{\it Splitness.---~} The main analogy with split reductive groups
is that the torus is a quotient of a finite numbers of copies of $\K^\times$ and each 
root group is isomorphic to the additive group of the groundfield $\K$. Combined with a
geometric consequence of the Moufang property, the latter property of root groups says: \pec 

{\bf Lemma.---~} \it 
Every building arising from a split Kac-Moody group is {\rm homogeneous}, that is thickness
is constant $($equal to $\# \K + 1)$ over all panels. 
\cqfd\rm\pec 

The buildings arising from a split Kac-Moody group will be referred to
as {\it (split) Kac-Moody buildings}. 

\gec 

\centerline{\bf 2.3. The adjoint representation}

\pec

Besides the geometry of buildings, an important tool to study Kac-Moody groups is a linear
representation which generalizes the adjoint representation of algebraic groups. Let us
fix a Kac-Moody root datum ${\cal D}$ and the corresponding Tits functor $\gkm$. 

\pec 

{\bf A. Formal sums and distribution algebras. } Let us go back to the situation of
2.1.C; we work in particular with the prenilpotent pair of roots $\{ a;b \}$. 
In the $\Z$-form $\ud$, we consider the subring 
${\cal U}_{[ a;b ]_{\hbox{\sevenrm lin}}}:= 
\ud \cap {\cal U}\g_{[ a;b ]_{\hbox{\sevenrm lin}}}$. 
Let us define the $R$-algebra $\widehat{\bigl( \uab \bigr)_R}$ of formal sums 
$\displaystyle \sum_{c \in \pn a + \pn b} r_c u_c$, with $u_c$ homogeneous of degree $c$ in
$\ud$ and $r_c \! \in \! R$. 

\pec 

{\bf Definition.---~} \it 
Denote by ${\goth U}_{[ a;b ]_{\hbox{\sevenrm lin}}}$ the group functor which to each ring
$R$ associates the subgroup $\langle \exp (re_c) \mid r \! \in \! R, c \! \in \!
\Delta^{re} \cap (\N a + \N b) \rangle < \widehat{\bigl( \uab \bigr)_R}^\times$. 
Here $\exp (re_c)$ is the formal exponential $\sum_{n \geq 0} r^n e_c^{[n]}$ and 
$\widehat{\bigl( \uab \bigr)_R}^\times$ is the multiplicative group of the $R$-algebra 
$\widehat{\bigl( \uab \bigr)_R}$. 
\rm\pec

There is a little abuse of notation when denoting by  
${\goth U}_{[ a;b ]_{\hbox{\sevenrm lin}}}$ a new object, but we will justify it soon. 
Thanks to a generalized Steinberg commutator formula, one can prove first: \pec 

{\bf Proposition.---~} \it 
The group functor ${\goth U}_{[ a;b ]_{\hbox{\sevenrm lin}}}$ is the functor of points of a
smooth connected group scheme of finite type over $\Z$, with Lie algebra 
$\g_{[ a;b ]_{\hbox{\sevenrm lin}}} \cap \ud$. Its value over $\C$ is the
complex algebraic group $U_{[ a;b ]_{\hbox{\sevenrm lin}}}$ and we have 
${\rm Dist }_e U_{[ a;b ]_{\hbox{\sevenrm lin}}} \cong 
{\cal U}\g_{[ a;b ]_{\hbox{\sevenrm lin}}}$. 
\cqfd\rm\pec

The notation ${\rm Dist }_e$ is for the {\it algebra of distributions supported at 
unity } of an algebraic group. This proposition is the place where a bit of algebraic
differential calculus [DG, II.4] comes into play. We will not detail the use of it,
but just say that it allows to prove that the group functors above and of proposition 2.1.C
coincide [R\'e2, 9.3.3]. In other words, this result leads to a concrete description of
the group functors occuring in the amalgam defining $\st$. 

\pec 

{\bf B. The adjoint representation as a functorial morphism. } For every ring $R$, the
extension of scalars from $\Z$ to $R$ of $\ud$ is denoted by $\bigl( \ud \bigr)_R$. Denote
by $\aut_{filt}(\ud)(R)$ the group of automorphisms of the $R$-algebra $\bigl( \ud
\bigr)_R$ preserving its filtration arising form that of $\ud$ and its ideal $\ud^+ \otimes_{\pz} R$. 
The naturality of scalar extension makes the assignement $R \mapsto  \aut_{filt}(\ud)(R)$ 
a group functor. We have [R\'e2, 9.5.3]: \pec

{\bf Theorem.---~} \it 
There is a morphism of group functors $\Ad: \gkm \to \aut_{filt}(\ud)$ characterized by: \pec 

\centerline{$\displaystyle \Ad_R \bigl( x_s (r) \bigr) = 
\sum_{n \geq 0} {(\ad e_s)^n \over n!} \otimes r^n, \quad 
\Ad_R \bigl( x_{-s} (r) \bigr) = 
\sum_{n \geq 0} {(\ad f_s)^n \over n!} \otimes r^n$,} 

\pec

\centerline{$\Ad_R \bigl({\cal T}_\Lambda(R) \bigr)$ \rm fixes $\bigl( {\cal U}_0 \bigr)_R$
\quad \rm and \quad $\Ad_R(h)(e_a \otimes r )= c_a(h) (e_a \otimes r)$,}

\pec

for $s$ in $S$, for each ring $R$, each $h$ in ${\cal T}_\Lambda(R)$ and each $r$ in $R$.
\rm\pec 

{\it Sketch of proof.---~} No use of group combinatorics is needed for this result.
Roughly speaking, to define the adjoint representation just consists in inserting \og
$\ad$\fg in the arguments of the formal exponentials for the root groups, and in making the
torus operate diagonally w.r.t. the duality given by ${\cal D}$. One has to follow all the
steps of the definition of Tits functors, and to verify that all the defining relations
involved are satisfied by the partial linear actions of the characterization. Proposition
2.3.B enables to justify that the first formula defines a representation of the Steinberg
functor. Then, the torus must be taken into account, but the verification concerning the
rest of the defining relations of $\gkm$ is purely computational. The naturality of the
representation is an easy point. 
\cqfd\pec

{\bf C. The adjoint representation over fields. } Restriction to values of Tits functors
over fields brings further properties, precisely because this allows the use of group
combinatorics according to 2.2.A -- see [R\'e2, 9.6] for proofs. 

\pec

{\bf Proposition.---~} \it 
{\rm (i) } Over fields, the kernel of the adjoint representation is the center of
the Kac-Moody group $G$. This center is the intersection of kernels of
characters of $T$, namely the centralizers in $T$ of the root groups. 

{\rm (ii) } If the groundfield is algebraically closed, the image of
$\Ad$ is again a Kac-Moody group. 
\cqfd\rm\gec 

\centerline{\bf 2.4. Algebraic subgroups}

\pec

One of the main interests of the adjoint representation is that it provides
algebro-geometric structures for the family of small subgroups of a Kac-Moody group. 

\pec 

{\bf A. Statement. } The following result is, up to technicalities, the combination of the
abstract Levi decompositions and restrictions of source and target of the adjoint
representation [R\'e2, 10.3]. 

\pec

{\bf Theorem (Algebraic structures).---~} \it
For every balanced subset $\Omega$, denote by $G(\Omega)$ the quotient group 
${\rm Fix}_G(\Omega)/Z \bigl( {\rm Fix}_G(\Omega) \bigr)$ of the
corresponding small subgroup by its center. Then $G(\Omega)$ admits a
natural structure of algebraic group, arising from the adjoint representation.
\cqfd\rm\pec 

{\it Remark.---~} An important question then is to know how far the quotient
$G(\Omega)$ is from the small subgroup ${\rm Fix}_G(\Omega)$ itself. To be more precise,
let us fix a pair $\A$ of opposite apartments containing the balanced subset $\Omega$.
Thanks to a suitable $g$ in $G$, we may assume $\A$ is the standard pair of opposite
apartments. At the group level, this corresponds to conjugating the small subgroup so as to
make it contain the standard Cartan subgroup $T$. Then the center of ${\rm Fix}_G(\Omega)$
is the intersection of the kernels (in $T$) of the characters associated to the roots in
$\Phi(\Omega)$. (These roots index the root groups generating the small subgroup together
with $T$). This shows that the quotient is always by a group abstractly isomorphic to a
diagonalisable algebraic group. This point will be important when we handle Cartan
subgroups, because preimages of diagonalisable groups will stay diagonalisable. 

\pec

{\it Examples.---~} Let us consider some examples from 1.5.C, i.e., from some
small subgroups of ${\rm SL}_2(\K[t,t^{-1}])$ and ${\rm SL}_3(\K[t,t^{-1}])$. 
There will be no quotient in the ${\rm SL}_2$ case: the
standard Cartan subgroup is one-dimensional. By taking a pair of
opposite middles of edges, one gets a small subgroup isomorphic to the
$M$-factor of the last example in 1.5.C. This group consists of the 
matrices $\pmatrix{X & 0 \cr 0 & {\rm det}X^{-1}}$, with $X \! \in \!
{\rm GL}_2(\fq)$. Its center is one-dimensional and this is the corresponding quotient that
actually admits an algebraic structure by the procedure above. 

\pec

{\bf B. } {\it Sketch of proof.---~} The main tools for this proof are basic results
from algebraic geometry or algebraic group theory, such as constructibility of images of
morphisms, use of tangent maps to justify that a bijective algebraic morphism is an
algebraic isomorphism...  Let us make the same reductions as in the remark above. In
particular, $\Omega$ is in the standard pair of opposite apartments.  Recall also the Levi
decomposition ${\rm Fix}_G(\Omega) = M(\Omega,\A) \ltimes U(\Omega,\A)$, with
$M(\Omega,\A)$ satisfying the abstract (TRD) axioms. 

\pec

The first step is the existence of a finite-dimensional subspace $W(\Omega)_{\pk} \subset
\bigl( \ud \bigr)_{\pk}$ such that the resctriction $\Ad_\Omega$ of the adjoint
representation is the center $Z \bigl( {\rm Fix}_G(\Omega) \bigr) (< T)$ 
[R\'e2, 10.3.1]. Consequently,
$\Ad_\Omega$ is faithful on $U(\Omega)$ and on the subgroup of $M(\Omega,\A)$ generated by a
positive half of $\Phi^m$. The latter group is the \og unipotent radical\fg of an
abstract  Borel subgroup of $M(\Omega,\A)$. Then, it is proved that the
isomorphic images of these two groups are unipotent closed $\K$-subgroups of ${\rm GL}
\bigl( W(\Omega)_{\pk} \bigr)$ [R\'e2, 10.3.2]. 

\pec 

The next step is to show that the abstract image 
$M_{\Omega,\hbox{\sevenbf A}}:= \Ad_\Omega \bigl( M(\Omega,\A) \bigr)$ is also a closed
subgroup.  This is done by proving that its Zariski closure in 
${\rm GL} \bigl( W(\Omega)_{\pk} \bigr)$ is reductive and comparing
the classical decompositions so obtained. The proof of this fact works as
follows. We know at this level that the image by $\Ad_\Omega$ of the abstract Borel subgroup
of $M(\Omega,\A)$ is closed. By Bruhat decomposition and local closedness of orbits, this
gives a bound on the dimension of the Zariski closure. The use of Cartier dual numbers
allows to compute the Lie algebras involved and to show that the bound is reached.  The
triviality of the unipotent radical follows from singular tori arguments [R\'e2,
10.3.3]. 

\pec

The final step is to show that the semi-direct product so obtained in 
${\rm GL} \bigl( W(\Omega)_{\pk} \bigr)$ is an algebraic one: this is proved by Lie
algebras (tangent maps) arguments [R\'e2, 10.3.4]. 
\cqfd\gec 

\centerline{\bf 2.5. Cartan subgroups}

\pec

We present now a general argument combining negative curvature and algebraic groups. 
Then we apply it to give a natural (more intrinsic) definition of  {\it Cartan subgroups},
previously defined by conjugation. Recall that a group of isometries is {\it bounded } 
if it acts with bounded orbits [Br2, p.160]. 

\pec

{\bf A.  A typical argument. } We suppose we are in the following setting: a
group $H$ acts in a compatible way on a Kac-Moody group $G$ and on its twin
building. (For instance, $H$ may be a subgroup of $G$ but this is not
the only example.) 

\pec 

{\bf Argument.---~} \it 
{\rm 1. } Justify that $H$ is a bounded group of isometries of both buildings
and apply the Bruhat-Tits fixed point theorem to get a balanced
$H$-fixed subset $\Omega$.

{\rm 2. } Apply results from algebraic groups to the 
small $H$-stable subgroup ${\rm Fix}_G(\Omega)$. 
\rm\pec

{\it Remark.---~} In subsection 3.2.B, the argument is applied to the 
Galois actions arising from the definition of almost split forms. 

\pec

{\bf B. Cartan subgroups. } Recall that the {\it standard Cartan
subgroup } is the value $T$ of the split torus ${\cal T}_{\cal D}$ on the grounfield $\K$.
The provisional definition of Cartan subgroups describes them as conjugates of $T$. We want
to give a more intrinsic characterization of them, by means of their behaviour w.r.t. the
adjoint representation $\Ad$. We suppose here that the groundfield $\K$ is infinite. 

\pec 

{\it A criterion for small subgroups.---~} This criterion precisely involves the
adjoint representation. Let $H$ be a subgroup of the Kac-Moody group $G$. Call $H$ {\it
$\Ad$-locally finite } if its image by $\Ad$ is locally finite, that is the linear
span of the $H$-orbit of each point in $\bigl( \ud \bigr)_{\pk}$ is finite-dimensional. 

\pec 

{\bf Theorem.---~} \it 
The following assertions are equivalent. 

\pec 

{\rm (i) } $H$ is $\Ad$-locally finite. 

{\rm (ii) } $H$ intersects a finite number of Bruhat double classes for each sign $\pm$. 

{\rm (iii) } $H$ is contained in the intersection of two fixators of spherical facets of
opposite signs. 
\rm\pec 

{\it Remark.---~} 1. The importance of this theorem lies in the fact that it is a
connection between the two $G$-actions we chose to study $G$: the
adjoint representation and the action on the twin building. 

2. This result is due to V. Kac and D. Peterson [KP3, Theorem 1] who proved a
lemma relating the $Q$-gradation and the Bruhat decompositions. This lemma was
used by them to prove conjugacy of Cartan subgroups over $\C$, with Kac-Moody groups
defined as automorphism groups of Lie algebras. No argument of negative curvature appears
in their work, they invoke instead convexity properties of a certain distance function on
the Tits cone. 

\pec 

{\it Definition and conjugacy of $\K$-split Cartan subgroups.---~} 
We say that a subgroup is {\it $\Ad$-diagonalizable } if its image by the
adjoint representation is a diagonalisable group of automorphisms of $\bigl( \ud
\bigr)_{\pk}$. Until the end of the section, we will assume the groundfield $\K$ to be 
infinite so as to characterize algebraic groups by their points. 

\pec 

{\bf Definition.---~} \it 
A subgroup of $G$ is a {\rm ($\K$-split) Cartan subgroup } if it is $\Ad$-diagonalizable 
and maximal for this property. 
\rm\pec 

We will omit the prefix \og $\K$-split\fg when the groundfield is separably closed.

\pec 

{\bf Theorem.---~} \it 
The $(\K$-split$)$ Cartan subgroups are all conjugates by $G$ of the standard Cartan subgroup
$T$. 
\rm\pec

{\it Sketch of proof. } [R\'e2, 10.4.2]---~Basically, this is the application of
argument 2.5.A, thanks to the smallness criterion above. An $\Ad$-diagonalisable subgroup
has to be small and its image in the algebraic subgroup is a Cartan subgroup in the algebraic
sense. We have then to use transitivity of $G$ on pairs of opposite apartments and the
algebraic conjugacy theorem of Cartan subgroups to send our subgroup in $T$ thanks to a
suitable element of
$G$. 
\cqfd\pec 

{\bf C. Connection with geometry. } Combining this theorem with proposition 2.2.B and
the fact that $G$ is transitive on the pairs of opposite apartments, we get [R\'e2,
10.4.3 and 10.4.4]: \pec

{\bf Corollary.---~} \it 
Suppose the groundfield $\K$ is infinite. Then, there is a natural $G$-equivariant
dictionary between the pairs of opposite apartments and the $\K$-split Cartan subgroups of
$G$. Moreover, given a facet $F$, this correspondence relates pairs of opposite apartments
containing $F$ and $\K$-split Cartan subgroups of $G$ contained in the parabolic subgroup 
${\rm Fix}_G(F)$. 
\cqfd\rm\pec

We saw in the abstract study of (TRD)-groups that parabolic subgroups do not form the
only interesting class of subgroups. Let us consider a balanced subset $\Omega$ of the twin
building of $G$. A natural question is to know what the correspondence says in this case. 
First we choose a pair of opposite apartments $\A$ containing
$\Omega$. We work in the
conical realization, that is see $\A$ as the union $\ct \cup -\ct$ of a Tits cone and its
opposite in the real vector space $V^*$. Then it makes sense to define the 
{\it convex hull } of $\Omega$ -- denoted by ${\rm conv}(\Omega)$ -- to be the trace on 
$\ct \cup -\ct$ of the convex hull determined in $V^*$. The same trick enables to define
the {\it vectorial extension } ${\rm vect}(\Omega,\A)$. The difference of notation is
justified by [R\'e2, 10.4.5]: \pec

{\bf Proposition.---~} \it 
{\rm (i) } $\K$-split Cartan subgroups of $G$ contained in ${\rm Fix}_G(\Omega)$ are the
preimages of the maximal $\K$-split tori of the algebraic group $G(\Omega)$. Besides, they
are in one-to-one $G$-equivariant correspondence with the pairs of opposite apartments
containing $\Omega$. 

{\rm (ii) } The subspace ${\rm conv}(\Omega)$ and the subgroup $U(\Omega)$ only depend on 
$\Omega$. 
\rm\pec

{\it Remark.---~} The independence of $\A$ for the convex hull is actually useful for
the Galois descent. 

\vfill\eject

\centerline{\bf 3. RELATIVE KAC-MOODY THEORY}

\gec 

We describe now an analogue of Borel-Tits theory. Section 3.1 defines the {\it almost split
forms } of Kac-Moody groups. This gives the class of groups which are concerned by the
theory.  In \S 3.2, we state the structure theorem for rational points and try to give an
idea of its quite long proof. The geometric method for the Galois
descent follows faithfully the lines of 
G. Rousseau's work, who used it for the characteristic 0 case [Rou1,2]. 
Section 3.3 introduces the class of {\it quasi-split } groups
and states that this is the only class of almost split Kac-Moody groups over finite fields. 
In \S 3.4, we focus on relative links and apartments, so as to compute thicknesses for
instance, and in \S 3.5 we use the classical example of the unitary group $SU_3$ to show how
things work concretely. 

\gec 

\centerline{\bf 3.1. Definition of forms}

\pec

The requirements are basically of two kinds. We have to generalize
conditions from algebraic geometry, and then to ask for isotropy
conditions. We suppose we are given a groundfield $\K$, we choose a separable closure $\ks$ 
(resp. an algebraic closure $\kb$) of it. We work with the group $G =
\gkm(\kb)$ determined by $\kb$ and the Kac-Moody root datum ${\cal
D}$. To ${\cal D}$ is also associated the $\Z$-form $\ud$; if $R$ is
any ring, then the scalar extension of $\ud$ from $\Z$ to $R$ will be denoted
by $\bigl( \ud \bigr)_R$. 

\pec 

{\bf A. Functorial forms. Algebraic forms. } Here is the part of the
conditions aiming to generalize algebraic ones. 
Unfortunately, $\K$-forms of Kac-Moody groups are not defined by a list of conditions
stated once and for all at the beginning of the study. 
One has to require a condition, then to derive some
properties giving rise to the notions involved in the next requirement, and so
on. We sum up here the basic steps of this procedure.

\pec 

{\it Use of functoriality.---~}  A Kac-Moody group is a value of a
group functor over a field. Consequently, it is natural to make occur a
group functor as main piece of a $\K$-form. A {\it functorial
$\K$-form } of $G$ is a group functor defined over field extensions of
$\K$ which coincides with the Tits functor $\ug_{\cal D}$ over extensions of
$\kb$. This is a really weak requirement since in the finite-dimensional
case, it just takes into account the functor of points of a scheme. 
Still, this condition enables to make the Galois group $\Gamma:={\rm
Gal}(\ks/\K)$ operate on $G$ [R\'e2, 11.1.2].
We explained in 2.2.A that the adjoint representation is a substitute for a global
algebro-geometric structure on a Kac-Moody group. It is quite natural then to work with
this representation to define forms [R\'e2, 11.1.3]. 

\pec 

{\bf Definition.---~} \it 
A {\rm prealgebraic $\K$-form } of $G$ consists in the datum of a
functorial $\K$-form $\ug$ and of a filtered $\K$-form $\uk$ of the
$\kb$-algebra $\ukb:= \ud \otimes_{\pz} \kb$, 
both satisfying the following conditions. 

\pec 

{\rm (PREALG1) } The adjoint representation $\Ad$ is
Galois-equivariant. 

{\rm (PREALG2) } The value of $\ug$ on a field extension gives rise to
a group embedding. 
\rm\pec 

{\it Remark.---~} It is understood in {\it filtered } that $\uk$
is a direct sum $\K \oplus \uk^+$, where $\uk^+$ is a $\K$-form of the
ideal $\bigl( \ud^+ \bigr)_{\pkb}$. The existence of a filtered $\K$-form of the
$\kb$-algebra $\ukb$ enables this time to define a $\Gamma$-action on
$\ukb$, and finally on the automorphism group $\autf(\kb)$. 

\pec 

{\it Example.---~} As a fundamental example, we can consider the
datum of the restriction of the Tits functor to the $\K$-extensions
(contained in $\kb$) and the $\K$-algebra $\bigl( \ud \bigr)_{\pk}$. 
This will be referred to as the {\it split form } of $G$. 

\pec

{\it Algebraicity conditions.---~} We suppose now we are given a
{\it splitting extension $\E/\K$ } of the prealgebraic form
$(\ug,\uk)$, that is a field over extensions of which the form is the
split one. We also assume that {\it the field $\E$ is infinite } and
that {\it the extension $\E/\K$ is normal}. The definition of Cartan
subgroups in terms of $\Ad$ and the $\Gamma$-equivariance of this
representation provides [R\'e2, 11.2.2]: \pec 

{\bf Lemma.---~} \it 
The image of an $\E$-split Cartan subgroup by a Galois automorphism is
still an $\E$-split Cartan subgroup. 
\cqfd\rm\pec 

In view of the conjugacy of Cartan subgroups, 
this lemma suggests to rectify Galois automorphisms so as to make them
stabilize the standard Cartan subgroup $T$. If $\gamma$ denotes a Galois
automorphism such that $\gamma T = g T g^{-1}$, then 
$\bar\gamma:={\rm int}g^{-1} \circ \gamma$ will denote such a rectified
automorphism. At the $\E$-algebras level, the symbol $\bar\gamma$ will
denote the rectified (semilinear) automorphism $\Ad g^{-1} \circ
\gamma$, where $\gamma$ is the Galois automorphism of $\bigl( \ud
\bigr)_{\pe}$ arising from the form. 
The choice of $g$ modulo $N$ will be harmless. Another class of interesting
subgroups is that of root groups, for which things are not so
nice. Still, we have [R\'e2, 11.2.3]: \pec

{\bf Lemma.---~} \it 
If the $c_s$'s are free in the lattice $\Lambda$ of the Kac-Moody root
datum ${\cal D}$, then the image of a root group $($w.r.t. to $T)$ under a rectified Galois
automorphism is still a root group $($w.r.t. to $T)$. 
\cqfd\rm\pec 

Note that the action on the group $G$ has no reason to respect its
nice combinatorics. The additional requirements are then [R\'e2, 11.2.5]: \pec 

{\bf Definition.---~} \it 
An {\rm algebraic $\K$-form } of $G$ is a prealgebraic form
$(\ug,\uk)$ such that for each $\gamma \! \in \! \Gamma$, one has: \pec 

{\rm (ALG0) } the rectified automorphism $\bar\gamma$ stabilizes 
the family of root groups $($w.r.t. to $T)$. 

{\rm (ALG1) } $\bar\gamma$ respects the abstract $Q$-gradation of 
$\bigl( \ud \bigr)_{\pe}$, inducing a permutation of $Q$ still denoted
by $\bar\gamma$ and satisfying the homogeneity condition 
$\bar\gamma(na) = n \bar\gamma (a)$, $a \! \in \! \Delta^{re}$, $n \!
\in \! \N$. 

{\rm (ALG2) } $\bar\gamma$ respects the algebraic characters and
cocharacters of $T$ among all the abstract ones. 
\rm\pec

{\it Remarks.---~} The algebraic characters (resp. cocharacters) are the ones arising
from $\Lambda$ (resp. $\tch{\Lambda}$): (ALG2) is indeed an algebraicity condition. 
Concerning (ALG1), the justification comes from algebraic differential calculus. If $G$
admitted an algebro-geometric structure -- a topological Hopf algebra of regular functions
for instance -- with $U_a$ as a closed subgroup, then the subring 
$\displaystyle \bigoplus_{j \leq n} \E {e_a^j \over j!}$ would be the $\E$-algebra of
invariant distributions tangent to $U_a$. So it is natural to require a rectified Galois
automorphism (stabilizing $T$) to respect the order of these distributions. Here are now the
main consequences of algebraicity [R\'e2, 11.2.5 and 11.3.2]. 

\pec 

{\bf Proposition.---~} \it 
{\rm (i) } $\bar\gamma$ is a group automorphism of $Q$, it stabilizes
$\Delta^{re}$ and induces the same permutation as the one defined via the root groups. 

{\rm (ii) } $\bar\gamma$ stabilizes $N$ inducing an automorphism of $W$ sending the
reflexion w.r.t. $\partial a$ to the reflection w.r.t. $\partial \bar\gamma(a)$. 

{\rm (iii) } If the Dynkin diagram of $A$ is connected, each Galois automorphism sends
a Borel subgroup on a Borel subgroup, possibly after opposition of its sign. 
\cqfd\rm\pec

In the non connected case, oppositions of signs can occur componentwise. At this step,
there is no reason why the $BN$-pair of fixed sign should be respected. This is a
requirement to be done. 

\pec 

{\bf B. Almost splitness. } This condition is the \og isotropy part\fg of the requirements
[Rou3]. 

\pec

{\bf Definition.---~} \it 
The algebraic form $(\ug,\uk)$ is {\rm almost split } if the Galois action stabilizes the
conjugacy of Borel subgroups of each sign. 
\rm\pec 

The stability of conjugacy classes of Borel subgroups enables to rectify each Galois
automorphisms $\gamma$ in such a way that $\gamma^*:={\rm int g_\sigma}^{-1} \circ \gamma$
stabilizes $B_+$ and $B_-$. This way, the element $g_\sigma$ is defined modulo 
$T = B_+ \cap B_-$. The assignement $\gamma \mapsto \gamma^*$ defines an $\Gamma$-action on
$W$ stabilizing $S$ and called the {\it $*$-action}. The terminology is justified by the
analogy with the action considered by J. Tits in the classification problem of semisimple
groups. Now we can turn for the first time to the other interesting $G$-space, namely the
twin building of the group [R\'e2, 11.3.3 and 11.3.4]. 

\pec 

{\bf Proposition.---~} \it 
{\rm (i) } Via group combinatorics, the Galois group operates on the buildings by
automorphisms, up to permutation of types by the $*$-action defined above. 

{\rm (ii) } The Galois group acts by isometries on the metric realizations of the building
of each sign. $\Gamma$-orbits are finite. 
\cqfd\rm\pec 

The last assertion is the starting point of the argument in 2.5.A: there exist
$\Gamma$-fixed balanced subsets. In fact, we can prove [R\'e2, 11.3.5]: \pec

{\bf Theorem.---~} \it 
Let us assume now $\E/\K$ Galois and let us consider $\Omega$, balanced and Galois-fixed.
Then the subspace $W(\Omega)_{\pe}$ is defined over $\K$, the homomorphism $\Ad_\Omega$ is
$\Gamma$-equivariant and its image is a closed $\K$-subgroup of 
${\rm GL} \bigl( W(\Omega)_{\pe} \bigr)$. 
\cqfd\rm\pec 

So to speak, this is the first step of the Galois descent: looking at
$\Gamma$-fixed balanced subsets and at the corresponding small subgroups
(resp. algebraic groups). 

\gec 

\centerline{\bf 3.2. The Galois descent theorem} 

\pec

This section states the main structure result about rational points of almost split
Kac-Moody $\K$-groups. The geometric side is at the same time the main tool and an
interesting result in its own right. 

\pec

{\bf A. Statement } [R\'e2, 12.4.3].  Here is the result which justifies the
analogy between the theory of almost split forms of Kac-Moody groups and Borel-Tits theory
for isotropic reductive groups. We keep the almost split $\K$-form $(\ug,\uk)$ of the
previous section.  We suppose now $\E=\ks$, i.e., that the form is split over the
separable closure
$\ks$. A Kac-Moody group endowed with such a form will be called an {\it almost split
Kac-Moody $\K$-group}. By definition, the {\it group of rational
points } of the form is the group of Galois-fixed points
$G(\K):=G^\Gamma$. 

\pec

{\bf Theorem (Galois descent). \trait} \it 
Let $G$ be a Kac-Moody group, almost split over $\K$. Then 

\pec

{\rm (i) } The set of $\Gamma$-fixed points in the twin building of $G$
is still a twin building ${\cal J}^\natural$. 
 
{\rm (ii) } The group of rational points $G(\K)$ is a $\trdl$-group 
for a natural choice of subgroups suggested by the
geometry of a pair of opposite apartments in ${\cal J}^\natural$. 
Indeed, ${\cal J}^\natural$ is a geometric realization of the twin building
abstractly associated to $G(\K)$. 

{\rm (iii) } Maximal $\K$-split tori are conjugate in $G(\K)$. 
\cqfd\rm\pec

{\it Remark.---~} The assumption of splitness over $\ks$ is not necessary in the
algebraic case, for which we can apply the combination of a theorem due to Cartier [Bo,
Theorem 18.7] and a theorem due to Chevalley-Rosenlicht-Grothendieck [Bo, Theorem 18.2]. 

\pec

{\bf B. } {\it Sketch of proof.---~} Let us start with notations
and terminology. The symbol $\Omega^\natural$ will always denote a
$\Gamma$-fixed balanced subset, and $\Omega$ will denote the
$\Gamma$-stable union of spherical facets covering it. 
We are working now in the conical realization of the twin building
so as to apply convexity arguments. 
A {\it generic subspace } is a subspace of a pair of opposite
apartments meeting a spherical facet. 

\pec 

{\it Relative geometric objects.---~} We define now relative objects with
suggestive terminology. This does not mean that we know at this step
that they form a building. 

\pec 

{\bf Definition.---~} \it
 {\rm (i) } A {\rm (spherical) $\K$-facet } is the subset of
$\Gamma$-fixed points of a $($spherical$)$ $\Gamma$-stable facet. 
A {\rm $\K$-chamber } is a spherical $\K$-facet of maximal closure. 
A {\rm $\K$-panel } is a $\K$-facet of codimension one in the closure
of a $\K$-chamber.   

{\rm (ii) } A {\rm $\K$-apartment } is a generic $\Gamma$-fixed
subspace, maximal for these properties. 
\rm\pec

{\it Remark.---~} These geometric definitions also give sense to
an opposition relation on $\K$-facets of oppposite signs, which will
be called {\it opposite} if they form a symmetric subset of a pair of
opposite apartments.  

\pec 

Let us turn now to the other application of argument 2.5.A,
with the Galois group as bounded isometry group [R\'e2, 12.2.1 and 12.2.3]. 

\pec 

{\bf Theorem.---~} \it 
For each $\Gamma$-fixed balanced subset $\Omega^\natural$, there
exists a $\Gamma$-stable pair of opposite apartments containing
$\Omega$. Since the Galois action preserves barycenters, 
${\rm conv}(\Omega^\natural)$ is $\Gamma$-fixed. 
\rm\pec

{\it Sketch of proof.---~} Theorem 3.1.B says that $\Ad_\Omega
\bigl( {\rm Fix}_G(\Omega) \bigr)$ is a $\K$-group, so 
it admits a Cartan subgroup defined over
$\K$. The $\Gamma$-equivariant dictionnary of 2.5.B says that this corresponds
geometrically to a $\Gamma$-stable pair of opposite apartments containing
$\Omega$. The last assertion follows from the fact that the convex
hull can be determined in any apartment containing $\Omega^\natural$. 
\cqfd\pec

The last assertion is the starting point of convexity arguments, which
enable to prove the following facts [R\'e2, 12.2.4]: \pec 

-- Two $\K$-facets of opposite signs are always contained in a
$\K$-apartment. 

-- There is an integer $d$ such that $\K$-apartments and $\K$-chambers are all
$d$-dimensional. 

-- The convex hull of two opposite $\K$-chambers is a
$\K$-apartment, and each $\K$-apartment can be constructed this way. In
particular, {\it $\K$-apartments are symmetric double cones}. 

-- Each $\K$-facet is in the closure of a spherical one. This enables to define $\K$-chambers as 
$\K$-facets ({\it a priori} not spherical) of maximal closure. 

\pec

{\it Rational transitivity results.---~} The next step is to study
${\cal J}^\natural$ as a $G(\K)$-set [R\'e2, 12.3.2]. 

\pec 

{\bf Theorem.---~} \it 
{\rm (i) } The group of rational points $G(\K)$ is transitive on the
$\K$-apartments. 

{\rm (ii) } Two $\K$-facets of arbitrary signs are always contained in
a $\K$-apartment. 
\cqfd\rm\pec 

{\it Remarks.---~} The proof of point (ii) requires a Levi
decomposition argument for a $\Gamma$-fixed balanced subset made of
three points, so Levi decompositions are useful for balanced subset
more general than pairs of spherical points of opposite signs. 

\pec 

{\it Relative roots. Anisotropic kernel.---~} 
To introduce a relative (TRD)-structure on the rational points, 
we need of course to define a relative Coxeter system
$(W^\natural,S^\natural)$. We can make the choice of 
a pair of $\K$-facets $\pm \fr$ of opposite signs and of a pair 
$\pm c$ of opposite chambers in such a way that: \pec

-- $\fr$ and $-\fr$ define a $\K$-apartment 
$\ak$ by $\ak:= {\rm conv}(\fr \cup -\fr)$; 

-- the $\Gamma$-stable facet $\pm F$ such that $\pm \fr = (\pm F)^\Gamma$
is in the closure of the chamber $\pm c$;

--  the pair of opposite apartments $\A$ defined by $\pm c$ is
$\Gamma$-stable. 

\pec 

We see $\A$ as the union of the Tits cone and its opposite in
the real vector space $V^*$. At last, we denote by $\lr$ the linear
span of $\ak$ in $V^*$. 

\pec 

{\bf Definition.---~} \it 
The restriction $\arel:= a \mid_{\lr}$ of a root $a$ of $\A$ $($seen as a linear
form on $V^*)$ is a {\rm $\K$-root } of $\ak$ if the trace $D(\arel):=\ak \cap D(a)$
is generic.  In this case, $D(\arel)$ is called the {\rm $\K$-halfspace } of
$\arel$. 
\rm\pec 

For a given $\K$-halfspace $D(\arel)$, we consider the set $\Phi_{\arel}$ of roots
$b$ with $D(b) \cap \ak = D(\arel)$. This set is stable under $\Gamma$ and so
is the group generated by the root groups $U_b$, $b \! \in \!
\Phi_{\arel}$. Let us denote by $V_{D(\arel)}$ the group of its
$\Gamma$-fixed points. 

\pec 

{\bf Definition.---~} \it 
{\rm (i) } The group $V_{D(\arel)}$ is the {\rm relative root group } associated to
the $\arel$. 

{\rm (ii) } The fixator  ${\rm Fix}_G(\ak)$ -- which we denote
by $Z(\ak)$ -- is the {\rm anisotropic kernel } associated to $\ak$. 
\rm\pec

{\it Rational group combinatorics.---~} 
The next step is to prove that $G(\K)$ admits a nice combinatorial
structure, stronger than $BN$-pairs but weaker than the
$\trdl$-axioms. 
This group combinatorics is due to V. Kac and D. Peterson: the {\it
refined Tits systems } [KP2]. This is an abstract way to obtain the relative
Coxeter system $(W^\natural,S^\natural)$ we are looking for, via
$BN$-pairs. We have [R\'e2, 12.4.1 and 12.4.2]: \pec

{\bf Lemma.---~} \it 
The group $W^\natural$ is a subquotient of $G(\K)$, namely the
quotient of the stabilizer in $G(\K)$ of $\ak$ by its fixator.
The half of each sign of the double cone $\ak$ is a geometric realization of the
Coxeter complex of $W^\natural$, and the latter group is 
generated by reflections w.r.t. to $\fr$. The $\K$-halfspaces $D(\arel)$
are in $W^\natural$-equivariant bijection
with the abstract set of roots $\Phi_{\pk}$ of $(W^\natural,S^\natural)$. 
\cqfd\rm\pec 

We have now an index set for the $\trdl$-axioms. What is still
unclear is how to define a linear interval of $\K$-halfspaces. For $D(\arel)$
and $D(\br)$, we set 

\pec 

\centerline{$[D(\arel);D(\br)]_{\hbox{\sevenrm lin.}}:= 
\{ D(c^\natural) \mid \exists \lambda, \mu \! \in \R_+ \quad c^\natural = 
\lambda\arel + \mu\br \}$.}

\pec 

With all the objects defined above, point (i) of the Galois descent theorem more
precisely says that {\it the group of rational points $G(\K)$ --
endowed with the $\K$-points  $Z(\K)^\Gamma$ 
of the anisotropic kernel and the family $\{ V_{D(\arel)} \}_{D(\arel) \in
\Phi_{\pk}}$ -- satisfies the $\trdl$-axioms}. Point (ii) says that
{\it the set of $\Gamma$-fixed points in the twin building is the 
realization of a twin building}. Besides, all the suggestive
terminology for relative objects is justified. 

\pec

{\it Maximal $\K$-split tori.---~} In the formulation above, we
forgot the maximal $\K$-split tori. According to 2.5.C, it is natural
to expect a relative dictionary between $\K$-apartments and maximal
$\K$-split tori. A {\it $\K$-split torus } of $G$ is a $\K$-split
subgroup of a $\Gamma$-stable Cartan subgroup (this definition makes
sense since the latter groups admit natural structures of $\K$-tori). 
Let us keep the $\K$-apartment $\ak$, and $Z(\K)={\rm Fix}_G(\ak)$. 

\pec 

{\bf Proposition.---~} \it 
The anisotropic kernel $Z(\ak)$ is the fixator of any pair of opposite 
$\K$-chambers $\pm \fr$ in $\ak$. It is a $\Gamma$-stable small
subgroup and the Levi factor of the $\Gamma$-stable parabolic
subgroup ${\rm Fix}_G(\pm \fr)$. Its associated algebraic group is
semisimple anisotropic. 
\rm\pec 

{\it Idea of proof. } [R\'e2, 12.3.2]---~The first assertion follows from example 2
in 1.5.B, and the anisotropy part is justified as in the proof of Lang's
theorem below (3.3.A). 
\cqfd\pec 

As the quotient of a small subgroup defining its associated algebraic
group is by a diagonalisable group, the anisotropy assertion gives
rise to a $G(\K)$-equivariant map 

\pec 

\centerline{$
\matrix{ \{ {\hbox{\rm $\K$-apartments}} \} & \longrightarrow & \{ {\hbox{\rm Maximal $\K$-split tori}} \} \cr 
\ak & \mapsto  & Z_d^0(\ak),}$} 

\pec 

where $Z_d^0(\ak)$ is the connected $\K$-split part of the center of $Z(\ak)$. A more
precise statement of point (iii) of the Galois descent theorem is that {\it each $\K$-split
torus can be conjugated by $G(\K)$ into $Z_d^0(\ak)$}. In particular, the map above is
surjective. In fact, it is bijective for $\K$ large enough. 

\gec

\centerline{\bf 3.3. A Lang type theorem}

\pec 

We describe here the specific situation in which the groundfield is finite. 

\pec

{\bf A. Quasi-splitness. } An algebraic $\K$-group is {\it
quasi-split } if it contains a Borel subgroup defined over $\K$, 
so in the twin situation, it is natural to introduce the following 

\pec

{\bf Definition.---~} \it 
A Kac-Moody $\K$-group is {\rm quasi-split } if it contains two opposite Borel subgroups
stable under the Galois group $\Gamma$. 
\rm\pec

{\it Remarks.---~} Geometrically, this means that two opposite chambers are
Galois-stable. 

\pec 

For comparison, the classical Lang theorem [Bo, 16.6] asserts that over a finite field,
an algebraic group has to be quasi-split. 

\pec 

{\bf Theorem (Lang's theorem for Kac-Moody groups).---~} \it 
Let $G$ be a Kac-Moody group, almost split over a finite field $\fq$. Then $G$ is in fact
quasi-split over $\fq$. 
\rm\pec

{\it Sketch of proof. } [R\'e2, 13.2.2] \trait Actually, the proof of this result makes
use of the algebraic one. The argument works as follows. Look at the Levi decomposition of a
minimal
$\Gamma$-stable parabolic subgroup $P$ of $G$. By minimality, the Levi factor is an
algebraic $\fq$-group which cannot contain a proper parabolic subgroup defined over $\fq$.
Indeed, it would otherwise lead to a $\Gamma$-stable parabolic subgroup smaller than $P$.
Consequently, the Levi factor -- which is reductive --  has no proper
parabolic subgroup defined over $\fq$, but admits an $\fq$-Borel subgroup: it must be a
torus. 
\cqfd\pec

{\bf B. Down-to-earth constructions of forms. } Restriction to quasi-split forms enables
to define concretely twisted Kac-Moody groups. In fact, groups of this kind were already
considered by J.-Y. H\'ee and J. Ramagge,  who obtained structure theorems for \og rational
points\fg without reference to or use of buildings [H1 and H3], [Ra1 and
Ra2]. These groups were studied in the wider context of generalized {\it Steinberg
torsions},  which allowed J.-Y. H\'ee to define analogues of Ree-Suzuki groups. 
If we stick to Galois torsion as above, we have the following
constructive result [R\'e2, 13.2.3]. 

\pec 

{\bf Proposition.---~} \it 
Consider a generalized Cartan matrix $A$ with Dynkin diagram $D$. 
Suppose we are given a Galois extension $\E/\K$ and a
morphism $*: {\rm Gal}(\E/\K) \to \aut(D)$. Then the simply
connected $($resp. adjoint$)$ split Kac-Moody group defined over $\E$
admits a quasisplit form over $\K$ with $*$-action the morphism above.   
\cqfd\rm\pec 

{\it Remark.---~} The Dynkin diagram of a generalized Cartan matrix is defined in
[K]. The result is also valid with a suitable notion
of automorphism for a Kac-Moody root datum. As we will see it as a
procedure to construct twisted geometries which only depends on $A$, we will not be
interested in these refinements. 

\gec 

\centerline{\bf 3.4. Relative apartments and relative links} 

\pec

We precise now the determination of the shape of the apartments and of the
relative links. 

\pec 

{\bf A. Geometry of relative apartments. } Our aim here is to derive 
some more information about the geometry of $\K$-apartments from the
$*$-action of the almost split form. This is to be combined with the
constructive procedure of quasi-split forms 3.3.B. 

\pec 

We are still working with the
datum of a pair of $\K$-facets $\pm \fr$ of opposite signs and of a pair 
$\pm c$ of opposite chambers s.t. $\fr$ and $-\fr$ define a $\K$-apartment 
$\ak$ by $\ak:= {\rm conv}(\fr \cup -\fr)$; the $\Gamma$-stable facet
$\pm F$ with $\pm \fr = (\pm F)^\Gamma$ is in the closure $\pm \overline c$ and 
the pair of opposite apartments $\A$ defined by $\pm c$ is $\Gamma$-stable. 
Recall that $\A$ is the union of the Tits cone and its opposite in
the real vector space $V^*$ and that $\lr$ is the linear span of $\ak$
in $V^*$. By transitivity of $G(\K)$, all $\K$-chambers are
fixed points in $\Gamma$-stable facets of fixed type $S_0$. 

\pec 

{\bf Lemma.---~} \it Denote by $\Gamma^*$ the
automorphism group of the standard chamber $c$ arising from the
$*$-action. Then $\lr$ admits the following definition by linear
equations.

\pec

\centerline{$\lr = \{ x \! \in \! V^* \mid a_s(x)=0 \ \forall s \! \in \! S_0 
\ \hbox{\rm and} \ a_s(x)=a_t(x) \ \hbox{\rm for} \
\Gamma^*s=\Gamma^*t \}$,}

\pec 

and the $\K$-chamber is a simplicial cone of dimension $d \geq
S^\natural$ defined by the same equations for $x$ in $c_+$. 
\cqfd\rm\pec 

This result [R\'e2, 12.6.1] enables to determine the shape of a $\K$-chamber or a
$\K$-apartment from the datum of the $*$-action. In other words, this is a concrete
procedure to determine a part of the geometry of the twisted twin building. To
get more precise information, we have to use finer theoretical results concerning
infinite Kac-Moody type root systems. Following the abstract results
in [Ba], the datum of $S_0$ and of the $*$-action determines a
{\it relative Kac-Moody matrix}, a relative version of the matrix $A$. Then, there is a
similar rule to deduce a Coxeter from it, and it happens to be naturally isomorphic to the
relative Weyl group $W^\natural$. So another piece of the twisted combinatorics
is theoretically computable. 

\pec

{\bf B. Geometry of relative residues. } Once we have determined the
geometry of the slices of the twisted twin building, one may want to know
more about the local geometry of it. This we will also accomplish thanks to
the knowledge of $S_0$ and of the $*$-action. In the split case, according to 2.2.B, the
fixator of a pair of opposite spherical facets $F_\pm$ is abstractly
isomorphic to a split reductive group. Besides, a realization of its twin
building is given by the union of the residues of $F_+$ and $F_-$. The
relative version of this is also valid. More precisely, let us
consider a pair of opposite spherical $\K$-facets $\pm \fr$ whose
union will be denoted by $\Omega^\natural$. It can be shown then that
the Galois actions on the $\Gamma$-stable residues yields a
natural $\K$-structure on the reductive algebraic groups to which the
small subgroup ${\rm Fix}_G(\Omega^\natural)$ is isomorphic. 
For instance, if we assume that $\pm \fr$ are  spherical $\K$-panels, 
this provides information about rational root groups. In fact, we
have [R\'e2, 12.5.2 and 12.5.4]: \pec 

{\bf Lemma.---~} \it 
{\rm (i) } The union of the sets of $\K$-chambers whose closure
contains $\pm \fr$ is a geometric realization of the twinning of an
algebraic reductive $\K$-group. 

{\rm (ii) } Rational root groups are isomorphic to root groups in
algebraic reductive $\K$-groups. In particular, the only possible 
proportionality factor $>1$ between $\K$-roots is $2$ and rational
root groups are metabelian. 
\cqfd\rm\pec 

Combined with the geometric consequence of the Moufang property, the
second point enables to compute thicknesses.

\gec

\centerline{\bf 3.5. A classical twin tree}

\pec

{\bf A. The split situation. } Recall that ${\bf SL}_3(\K[t,t^{-1}])$
is a split Kac-Moody group defined over the field $\K$ and of affine type
$\widetilde A_2$. The Dynkin diagram is the simply-laced triangle. 
The simple roots of the $\widetilde A_2$-affine root system are
denoted by $a_0$, $a_1$ and $a_2$, for which we choose respectively the parametrizations 
$x_0 \!: r \mapsto \pmatrix{1&0&0\cr 0&1&0\cr -rt&0&1}$, 
$x_1 \!: r \mapsto \pmatrix{1&r&0\cr 0&1&0\cr 0&0&1}$ and 
$x_2 \!: r \mapsto \pmatrix{1&0&0\cr 0&1&r\cr 0&0&1}$. 
The standard Cartan subgroup is $T = \{ D_{u,v}:=
\pmatrix{u&0&0\cr 0&u^{-1}v&0\cr 0&0&v^{-1}} \mid u,v \! \in \! \E^\times \}$. 
At last, we work for this example with the Bruhat-Tits realization of the buildings. 

\pec 

{\bf B. The involution. } Let us consider a separable quadratic extension $\E/\K$.
We want to define a quasi-split $\K$-form of ${\rm SL}_3(\E[t,t^{-1}])$ so as to
produce a semi-homogeneous twin tree by quasi-split Galois descent. 
Fom a theoretical point of view (3.3.B), 
it is enough to find an involution of the Dynkin diagram. We
choose the inversion of types 1 and 2. We describe now concretely to
which Galois action this choice leads. Let $\iota$ (resp. $\tau$) be the
inversion (resp. the \og antidiagonal\fg symmetry) of matrices. 

\pec 

\centerline{$\tau: \pmatrix{a&b&c\cr d&e&f \cr g&h&i} \mapsto 
\pmatrix{i&f&c\cr h&e&b \cr g&d&a}$.}

\pec

Let $\cdot^\sigma$ be at the same time the nontrivial element of 
${\rm Gal}(\E/\K)$ and this operation on all matrix coefficients. 
Then the involution $\cdot^*:=\tau \circ \iota \circ \cdot^\sigma$ of 
${\bf SL}_3(\E[t,t^{-1}])$ is the one given by the quasi-split $\K$-form. 
It is completely determined by the formulas 
$x_1(r)^*=x_2(r^\sigma)$, $x_2(r)^*=x_1(r^\sigma)$ and $x_0(r)^*=x_0(r^\sigma)$,  
the analogues for negative simple roots, and $D_{u,v}^*=D_{v^\sigma,u^\sigma}$. The group
of monomial matrices lifting modulo $T$ the affine Weyl group $W={\cal S}_3 \ltimes \Z^2$ 
is also $\cdot^*$-stable. This implies the $\cdot^*$-stability of the standard apartments 
$\A_+=\{ wB_+ \}_{w \in W}$ and $\A_-=\{ wB_- \}_{w \in W}$. 
The set of $^*$-fixed points inside is a $\K$-apartment, since according
to the permutation of types, there cannot be a two-dimensional
$\K$-chamber. This set $\ak$ is the median of the standard equilateral triangle 
(intersecting the edge of type $0$) w.r.t. to which $\cdot^*$ operates
by symmetry in $\A$.  Consequently, each building is a tree and the relative Weyl group $W^\natural$
is the infinite  dihedral group $D_\infty$. Translations are lifted by the matrices 
$\pmatrix{ at^n &0&0 \cr 0&1&0 \cr 0&0& (at^n)^{-1} }$; reflexions by 
the matrices $\pmatrix{0&0& (at^n)^{-1} \cr 0&1&0 \cr -at^n &0&0}, 
(n \! \in \Z, \ a \! \in \! \K^\times)$. The reason why we call the
twin building so obtained a classical twin tree is that, up to sign
conventions, the Steinberg torsion we defined is that of the unitary group ${\rm SU}_3$
considered in [BrT2, 4.1.9]. 

\pec 

{\bf C. Semi-homogeneity. } What is left to do is to compute valencies. 
We will get semi-homogeneous trees since there are two shapes of intersections of the  
$\K$-apartment with walls. We sketch arguments for positive standard vertices
and edges, but the situation is completely symmetric w.r.t. signs and
can be deduced for every facet by transitivity of the rational points. 

\pec

The first (simplest) situation is given by the intersection of the standard
edge of type 0 with $\ak$. In this case, only one wall passes through
the intersection point and each of the two opposite $\K$-roots is the
trace on $\ak$ of a single root. The corresponding root groups are
stable and their groups of rational points are isomorphic to the
additive group of $\K$. Hence a valency at all transforms of this
$\K$-vertex equal to $1+\#\K$. 

\pec

The other situation corresponds to the fact that $\ak$ contains the
standard vertex of type $\{ 1;2 \}$ for which 3 roots leave the same $\K$-root
as trace on $\ak$. The small subgroup fixing this vertex and its
opposite in $\A$ is stable and abstractly isomorphic to ${\rm SL}_3$ over $\E$. 
According to the action around the vertex,
the three positive root groups passing through it are permuted and the
small subgroup is abstractly isomorphic to
a quasi-split reductive $\K$-group of rank one. 
The root groups of such groups are described in detail in [BrT2,
4.1.9]; they are rational points of the 
unipotent radical of a Borel subgroup of $SU_3$, in bijection
with $\K^3$. This can also be seen directly on the explicit
description of the Galois action above. Hence a valency at all transforms of this
$\K$-vertex equal to $1+\#\K^3$. 

\pec

This shows that we built a semi-homogeneous twin tree of valencies
$1+\#\K$ and $1+\#\K^3$. So in the case of a finite groundfield $\K=\fq$, one gets a
locally finite semi-homogeneous twin tree of valencies $1+q$ and $1+q^3$. 

\pec 

{\it Figure.---~} [Here $q=2$]

\centerline{\epsfysize=50mm $$\epsfbox{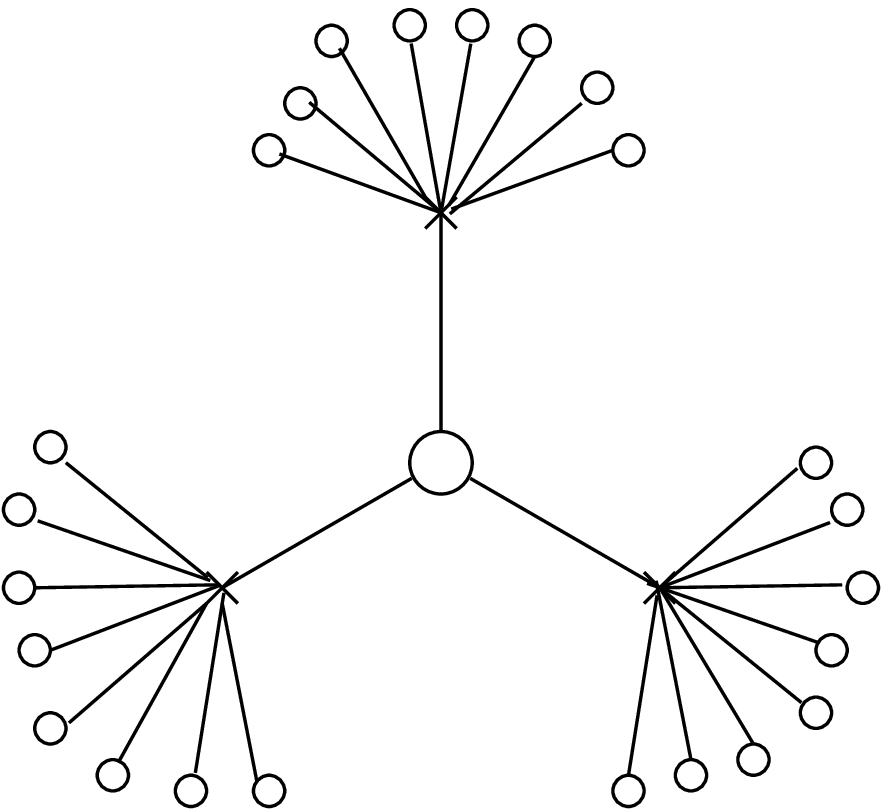}$$}

\vfill\eject

\centerline{\bf 4. HYPERBOLIC EXAMPLES. LATTICES}

\gec

So far, we just considered examples of Chevalley groups over Laurent polynomial rings
or twisted versions of them. Section 4.1 deals with split and twisted Kac-Moody groups whose
buildings are hyperbolic.  In \S 4.2, we show that Kac-Moody theory provides analogues of
arithmetic groups over function fields. 

\gec

\centerline{\bf 4.1. Hyperbolic examples}

\pec 

At last, some more exotic Kac-Moody geometries will be described. 

\pec

{\bf A. Hyperbolic buildings. } The definition of such spaces [GP] is the
metric version of that of buildings from the apartment system point of view. 

\pec 

{\bf Definition.---~} \it
Let $R$ be a hyperbolic polyhedron providing a Poincar\'e tiling of the hyperbolic
space ${\Bbb H}^n$. A {\rm hyperbolic building of type $R$ } is a piecewise polyhedral
cell complex, covered by a family of subcomplexes -- the {\rm apartments}, all isomorphic
to the tiling above and satisfying the following incidence axioms. 

\pec

{\rm (i)} Two points are always contained in an apartment. 

{\rm (ii)} Two apartments are isomorphic by a polyhedral arrow fixing their
intersection. 

\pec

A building whose apartments are tilings of the hyperbolic plane ${\Bbb H}^2$ will be
called {\rm Fuchsian}. 
\rm\pec

The main point about hyperbolic buildings is that we get this way ${\rm CAT}(-1)$-spaces 
i.e., metric spaces where geodesic triangles are at least as thin as in ${\Bbb H}^2$.
In particular, they are hyperbolic in the sense of Gromov. The additional incidence axioms
make these spaces nice geometries where to refine results known for hyperbolic metric
spaces, or to generalize results about symmetric spaces with (strict) negative curvature. 

\pec

{\it Remark.---~} Note that in contrast with the affine case, for hyperbolic buildings 
there is usually no more relation between the dimension of the apartments and the rank of
the Weyl group of the building. The two-dimensional case with a regular right-angled $r$-gon
is particularly striking: $r \geq 5$ may be arbitrary large. Nevertheless,
we will sometimes be interested in the case called {\it compact hyperbolic} by Bourbaki. It
corresponds to a tesselation of ${\Bbb H}^n$ by a hyperbolic simplex. In this particular
case, the relation alluded to above is valid; and the Weyl group shares another property
with the affine case: every proper subset of canonical generators gives rise to a finite
Coxeter subgroup. 

\pec

{\bf B. Split twin hyperbolic buildings. } Our purpose here is to show which kind of
hyperbolic buildings can be produced via Kac-Moody machinery. Let us make preliminary
remarks, valid in both split and twisted cases. First, the rule of point (i) in the
second definition of 2.1.A shows that the Coxeter coefficients involved in Kac-Moody theory
are very specific: 2, 3, 4, 6 or $\infty$. So by far not any Coxeter group appears as
the Weyl group of a generalized Cartan matrix. Besides, the thicknesses will always be of
the form \og 1 + prime power\fg, and will be constant in the split case. These are strong
conditions, but not complete obstructions. 

\pec 

{\it Hyperbolic realizations.---~} It is clear that standard gluing techniques 
[D2, \S 10] apply in the hyperbolic case once affine tilings replaced by hyperbolic
ones. Given $P$ a hyperbolic Poincar\'e polyhedron with dihedral angles $\neq 0$ are 
of the form $\pi/m$ with $m=2$, $3$, $4$ or $6$, and which may have
vertices at infinity; and given $p^k$ a prime power, 
there exists a Kac-Moody building of constant thickness
$1+p^k$ admitting a metric realization whose apartments are Poincar\'e tilings defined by
$P$ [R\'e1 \S2]. This building is a complete geodesic proper ${\rm CAT}(-1)$ metric
space.  This is just an existence result, it is expected to be combined with a uniqueness
assertion, but the following result due to D. Gaboriau and F. Paulin
([GP], th\'eor\`eme 3.5) is quite disheartening. 

\pec 

{\bf Theorem.---~} \it
Let $P$ be a regular hyperbolic Poincar\'e polygon with an even number
of edges and angles equal to $\pi/m, m \geq 3$. 
Let $L$ be a fixed generalized $m$-gon of classical type and sufficiently large characteristic. 
Then, there exist uncountably many hyperbolic buildings with
apartments the tiling provided by $P$ and links all isomorphic to $L$.
\cqfd\rm\pec

The latter theorem is in contrast with the following characterization proved
by M. Bourdon [Bou1, 2.2.1], which shows that in the two-dimensional case, there is 
some kind of  dichotomy according to whether $m \geq 3$ or $m=2$. The
latter kind of buildings will be referred to as the class of {\it right-angled
Fuchsian buildings}. 

\pec

{\bf Proposition.---~} \it 
Consider the Poincar\'e tiling of the hyperbolic plane ${\Bbb H}^2$ by
a right-angled $r$-gon $(r \geq 5)$. Then, there exists
a unique building $I_{r,q+1}$ satisfying the following conditions: \pec

\pech -- apartments are all isomorphic to this tesselation;

\pech -- the link at each vertex is the complete bipartite graph of parameters
$(q+1,q+1)$. 
\cqfd\rm\pec

Consequently, we know by uniqueness that the building $I_{r,q+1}$ -- $r \geq 5$, $q$  
a prime power -- comes from a Kac-Moody group. 

\pec 

{\bf C. Quasi-split Kac-Moody groups and hyperbolic buildings. } On the one hand, according
to subsection 2.2.B, there is no possibility to construct buildings with non-constant
thickness  via split Kac-Moody groups. On the other hand, we know by 3.5 that we can obtain
semi-homogeneous twin trees from classical groups. Here is an example of what can be
done in the context of hyperbolic buildings. 

\pec

{\it Semi-homogeneous trees in Fuchsian buildings. } [R\'e2, 13.3]---~Definition of $\fq$-forms relies on the existence of symmetries in the
Dynkin diagram. Let us sketch the application of method 3.3.B to the buildings
$I_{r,q+1}$.  Denote by $D$ the Coxeter diagram of the reflection group associated to the
regular right-angled $r$-gon $R \subset {\Bbb H}^2$. To get it, draw the $r$-th
roots of 1 on the unit circle and connect each pair of non consecutive points by an edge 
$\displaystyle {\buildrel \infty \over {\hbox{--------}}}$. To keep all symmetries, lift
this Coxeter diagram into a Dynkin diagram replacing the edges 
$\displaystyle {\buildrel \infty \over {\hbox{--------}}}$ by $\Longleftrightarrow$. 
The corresponding generalized Cartan matrix $A$ is indexed by $\Z/r$; it is defined by 
$A_{\bar i, \bar i}=2$, $A_{\bar i, \bar i+1}=0$ and -2's elsewhere. 

\pec 

{\it Figure.---~} [Here $r=6$] 

\centerline{\epsfysize=30mm $$\epsfbox{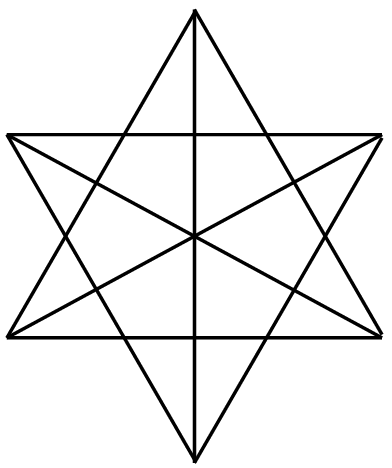}$$}

We suppose now that the integer $r$ is odd, so that $D$ is symmetric w.r.t. each axis
passing through the middle of two neighbors and the opposite vertex. We fix such a
symmetry, and consider the quadratic extension $\fqd/\fq$: this defines a quasi-split
form as in 3.5.B. Moreover, the arguments to determine the shape of an
$\fq$-chamber or an $\fq$-apartment are the same. We get a straight line cut into segments
and acted upon by the infinite dihedral group $D_\infty$. The twisted buildings are trees. 

\pec 

{\it Figure.---~} [Here $r=5$ and $q=2$]

\centerline{\epsfysize=25mm $$\epsfbox{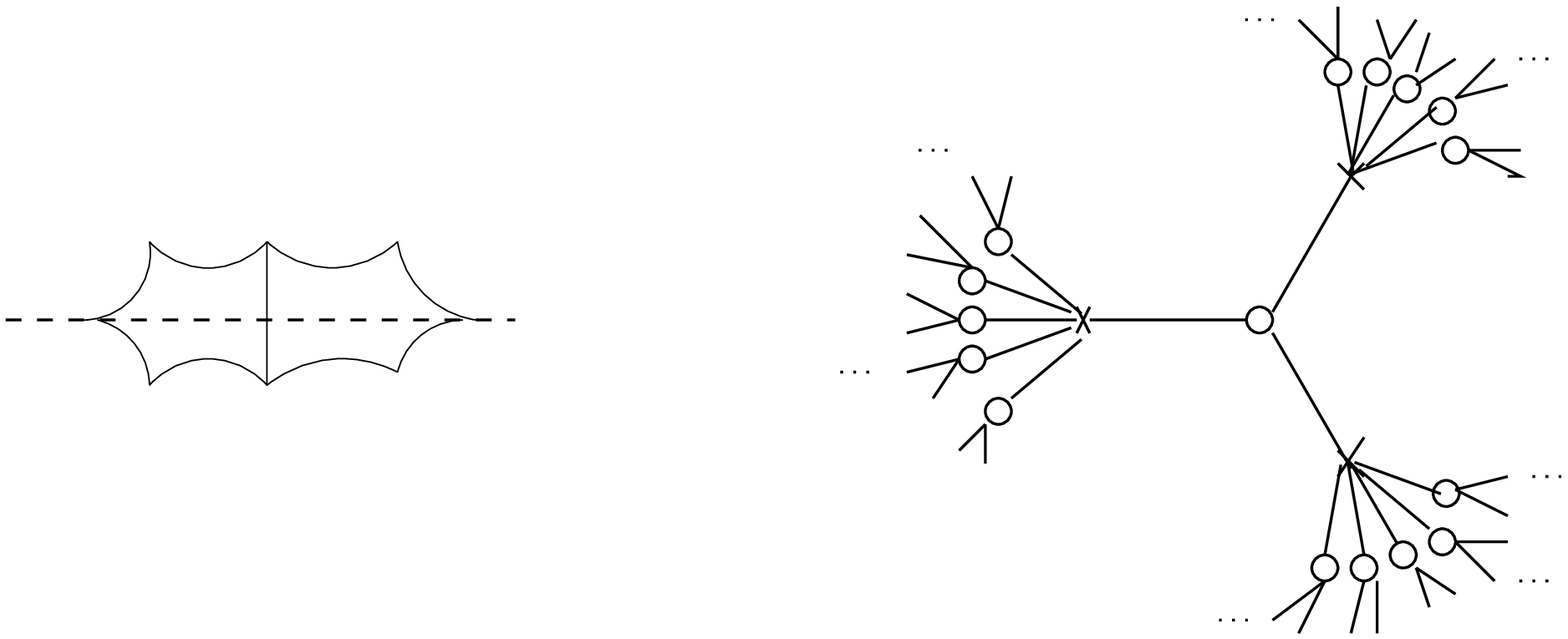}$$}

For the valencies, we must once again distinguish two kinds of intersection of the boundary
of
$R$ with the $\fq$-apartment, which leads to two different computations. When the
$\fq$-apartment cuts an edge, the rational points of the corresponding root group
are isomorphic to $\fq$, hence a valency equal to $1+q$. The other case is when the 
$\fq$-apartment passes through a vertex contained in two orthogonal walls, with two roots
leaving the same trace on the Galois-fixed line. The corresponding Galois-stable group is
then a direct product of two root groups (each isomorphic to $\fqd$) acted upon by interversion
and Frobenius: the valency at these points is then $1+q^2$. Consequently, we obtain this way
a semi-homogeneous twin tree with valencies $1+q$ and $1+q^2$. 

\pec

{\it Twisted Fuchsian buildings in higher-dimensional hyperbolic buildings. } 
[R\'e1, \S 4.5]---~Galois descent makes the dimension of apartments decrease, so if we want Fuchsian buildings
with non constant thickness, we must consider higher dimensional hyperbolic buildings. 
Let us consider a specific example obtained from a polyhedron in the hyperbolic ball 
${\bf B}^3$ and admitting many symmetries: it is the truncated regular hyperideal octahedron
with dihedral angles $\pi/2$ and $\pi/3$. This polyhedron $P$ is obtained by moving at
same speed from the origin eight points on the axis of the standard orthogonal basis of
$\R^3$. The condition on the angles impose to move the points outside of ${\bf B}^3$. What
is left to do then is to cut the cusps by totally geodesic hyperplanes centered at the
intersections of ${\bf S}^2=\partial{\bf B}^3$ with the axis. A figure of this polyhedron
is available in [GP]. The polyhedron $P$ admits a symmetry w.r.t. hyperplanes
containing edges and cutting $P$ along an octogon. The residue at such an edge -- where the
dihedral angle is $\pi/3$ -- is the building of ${\rm SL}_3$. This enables to define a
quasi-split $\fq$-form whose building is Fuchsian. The $\fq$-chambers are octogons and
thicknesses at $\fq$-edges are alternatively $1+q$ and $1+q^3$. 

\pec 

{\bf D. Hyperbolic buildings from a geometrical point of view. } Here
is now a small digression concerning single hyperbolic buildings, without the
assumption that they arise from Kac-Moody theory. This often implies
to restrict oneself to specific classes of hyperbolic buildings, but
provides stimulating results. 

\pec 

{\it Existence, uniqueness and homogeneity.---~} We saw in the
previous subsection that unicity results were useful to know which
buildings come from Kac-Moody groups. Nevertheless, the \og disheartening theorem\fg above
shows that the shape of the chambers and the links are not sufficient in this direction. F.
Haglund [Ha] determined conditions to be added so as to obtain characterization results
for two-dimensional hyperbolic buildings. His study enables him to know when the
(type-preserving) full automorphism group is transitive on chambers.  

\pec 

{\it Automorphism groups.---~} The first step towards group
theory then is to consider groups acting on a given hyperbolic
building. Two kinds of groups are particularly interesting: full automorphism groups and lattices. Concerning full
automorphism groups, F. Haglund and F. Paulin [HP] proved 
abstract simplicity and non-linearity results for two-dimensional
buildings. Besides, \'E. Lebeau [L] proved that some others of these groups are
strongly transitive on their buildings, hence admit a natural
structure of $BN$-pair -- see [Ron], \S 5. Concerning lattices, the
first point is the existence problem. In fact, the main construction technique
is a generalization in higher dimension of Serre's theory of graph of
groups [Se]. D. Gaboriau and F. Paulin [GP] determined the precise local
conditions on a {\it complex of groups } in the sense of A. Haefliger
[BH] to have a hyperbolic building as universal covering. This
method has the advantage to produce at the same time a uniform
lattice with the building. The lattices we will obtain by Kac-Moody
theory are {\it not } uniform. 

\pec 

{\it Boundaries: quasiconformal geometry and analysis.---~} 
The reason why boundaries of metric hyperbolic metric spaces are
interesting is that they can be endowed with a quasiconformal
structure -- see [GH], \S 7. In the case of hyperbolic buildings, the
structure can be expected to be even richer or better known, in view
of the ${\rm CAT}(-1)$-property and the incidence axioms. At least
for Fuchsian buildings, this is indeed the case as shown by works due
to M. Bourdon. In this situation, he relates the precise conformal
dimension of the boundary and the growth rate of the Weyl group [Bou2]. In a joint paper
with H. Pajot [BP1], the existence of Poincar\'e inequalities and the structure
of Loewner space is proved for the boundaries of the $I_{r,q+1}$'s. 

\pec 

{\it Rigidity.---~} Rigidity is a vast topic in geometric group
theory. In this direction, results were proved by M. Bourdon
first [Bou1] and then by M. Bourdon and H. Pajot [BP2]. The
first paper deals with Mostow type rigidity of uniform lattices. It
concerns right-angled Fuchsian buildings. A basic tool there is a
combinatorial crossratio constructed thanks to the building
structure involved. It enables to generalize arguments of hyperbolic 
geometry by D. Sullivan to the singular case of (boundaries of) buildings. 
The second paper uses techniques from analysis on singular spaces to prove a stronger
result implying Mostow rigidity, and in the spirit of results by P. Pansu [P1], 
B. Kleiner and B. Leeb [KL] for Riemannian symmetric spaces
(minus the real and complex hyperbolic spaces) and higher rank Euclidean buildings.
More precisely, M. Bourdon and H. Pajot proved [BP2]: \pec 

{\bf Theorem (Rigidity of quasi-isometries).---~} \it 
Every quasi-isometry of a right-angled Fuchsian building lies within
bounded distance from an isometry. 
\cqfd\rm\pec 

For the notion of a {\it quasi-isometry}, see for instance
[GH, \S 5]. Roughly speaking, it is an isometry up to an additive and a
multiplicative constant, surjective up to finite distance. 

\gec

\centerline{\bf 4.2. Analogy with arithmetic lattices}

\pec 

The purpose of this section is to present several arguments to see
Kac-Moody groups over finite fields as generalizations of arithmetic
groups over function fields. 

\pec

{\bf A. Kac-Moody topological groups and lattices. } Let $G$ be a Kac-Moody group, almost
split -- hence, quasi-split (3.3.A) -- over the finite field $\fq$. 
We denote by $\Delta_\pm$ the locally finite metric realizations of 
the buildings associated to $G$. Recall that all points there are of
spherical type, hence the fixator of any point is a spherical parabolic subgroup. 
By local finiteness, the (isomorphic) corresponding full automorphism groups are locally compact.
We use the notations $\aut_\pm:= \aut(\Delta_\pm)$. The groups $\aut_+$ and $\aut_-$ are
besides unimodular and totally disconnected. We choose a left
invariant Haar measure ${\rm Vol}_\pm$ and normalize it by 
${\rm Vol}_\pm \bigl( {\rm Fix}_{\paut_\pm}(c_\pm) \bigr) = 1$, 
where $c_\pm$ is the standard chamber.  

\pec 

Let us turn now to the analogy with arithmetic groups [Mar, p.61]. The affine case is
indeed a guideline since Chevalley groups over Laurent polynomials are 
$\{ 0;\infty \}$-arithmetic groups over the global field $\fq(t)$. 
This is the example of an arithmetic group for two places. If we
consider now the fixator of a facet, we get a group commensurable with
the value of the latter Chevalley group over the polynomials $\fq[t]$, so that each
fixator of facet is an arithmetic group for one place. 
Besides, the normalization we chose for the invariant measure
consists in the affine case in assigning volume 1 to an
Iwahori subgroup. We just replaced here a $p$-adic Lie group by the full automorphism group
of a building. 

\pec 

If we consider now the general case, here is a result supporting the analogy, 
independently proved by L. Carbone and H. Garland in the split case [CG, R\'e3]. 

\pec 

{\bf Theorem.---~} \it 
Let $G$ be an almost split Kac-Moody group over a finite field
$\fq$. Denote by $\displaystyle \sum_{n \geq 0} d_n t^n$ the 
growth series of the Weyl group $W$ of $G$. 
Assume $\displaystyle \sum_{n \geq 0}{d_n\over q^n} < \infty$. Then: 

\pec

{\rm (i)} The group $G$ is a lattice of (the full automorphism group of) 
$\Delta_+ \times \Delta_-$. 

{\rm (ii)} For any point in $x_- \! \in \! \Delta_-$, ${\rm Fix}_G(x_-)$ is 
a lattice of $\Delta_+$. 

{\rm (iii)} These lattices are not uniform. 
\cqfd\rm\pec

{\it Remarks.---~} 1. Non-uniformness is due to J. Tits, who
determined fundamental domains for actions of parabolic subgroups in a combinatorial
setting more general than the (TRD) axioms [T3, A1 \S 3]. 

2. Of course, the situation is completely symmetric w.r.t. signs: point (ii) is still
true after inversion of $+$ and $-$.  

3. There is a straightforward generalization of this result to the
case of locally finite Moufang twin buildings. The cardinal of the groundfield has to
be replaced by the maximal cardinal $q$ of root groups. According to the Moufang property
(1.3.B), this integer is geometrically computed as follows. Consider any chamber, then
$q+1$ is the maximal thickness among all the panels in its closure. 

\pec 

{\bf B. } {\it Sketch of proof.---~} We just consider the case of the fixator of a negative
point -- spherical negative parabolic subgroup. As we work over a finite field, we can make
use of the following facts due to Levi decompositions, 1.5.B and 1.5.A respectively. 

\pec 

(i) The fixator of a pair of (spherical) points of opposite signs is finite. 

(ii) Borel subgroups are of finite index in the spherical parabolic subgroups containing
them. 

\pec 

Discreteness follows from point (i), and according to point (ii), it is enough to show that
the fixator $B_-$ of a negative chamber is of finite covolume. The idea is to adapt Serre's
criterion on trees to compute a volume. This requires to look at the series 
$\displaystyle \sum_{d_+ \in {\cal C}} {1 \over \# {\rm Stab}_{B_-} (d_+)}$, where 
${\cal C}$ is a complete set of representatives of chambers modulo the discrete group
considered. According to Tits' result alluded to
above, one has ${\cal C}=\A_+$. Thanks to a special case of the Levi
decompositions (1.5.B),
the series is $\displaystyle \leq {1 \over \# T} \sum_{w \in W} {1 \over q^{\ell(w)}}$. 
Rearranging elements of $W$ w.r.t. length in the dominating series, one gets the series appearing in the
statement. 
\cqfd\pec

{\it Remark.---~} When $G$ is split over $\fq$, the assumption on the numerical
series is sharp. 

\gec

{\bf C. Classical properties for Kac-Moody lattices. } The result above should be seen as a
starting point. Indeed, one may wonder whether some classical properties of arithmetic
groups are relevant or true for this kind of lattices. Some answers are already known. 

\pec

{\it Moufang property.---~} 
The first remark in order for the analogy is to consider Moufang property. Recall that this
property is already of interest in classical situations but is not satisfied by all
$p$-adic Lie groups. For instance, it can be proved for many Bruhat-Tits buildings arising
from reductive groups over local fields of equal (positive) characteristic, but it is
known for no local field of characteristic zero [RT1, \S 2]. This remark suggests that
the analogy with arithmetic groups in characteristic 0 will not be fruitful. 

\pec

{\it Kazhdan's property {\rm (T)}.---~} 
This property has proved to be useful in many situations, and makes
sense for arbitrary topological groups, endowed with the discrete
topology for instance. In the classical case of Lie groups and their lattices, 
its validity is basically decided by the rank of the ambient Lie
group [HV, 2.8]. In rank $\geq 2$, all Lie groups are (T) and closed subgroups
of finite covolume inherit the property [HV, 3.4]. Over non-archimidean local
fields and in rank 1, no group is (T) since a necessary condition is to
admit a fixed point for each isometric action on a tree [HV, 6.4]: contradiction
with the natural action on the Bruhat-Tits tree. Another necessary
condition is for instance compact generation, that is finite
generation in the discrete case. 

\pec 

Property (T) is often  characterized in terms of unitary
representations, but the former necessary condition involving actions on
trees suggests to look for geometrical criteria having the
advantadge to make generalizations easier. An important example of this approach works as
follows.   On the one hand, H. Garland used the combinatorial Laplace operator on
Bruhat-Tits buildings  to obtain vanishing results for cohomology with coefficients in
finite-dimensional representations [G]. On the other hand, property (T) is
characterized by cohomology vanishing for every unitary
representation. The idea to generalize Garland's method to
infinite-dimensional representations is due to Pansu-Zuk [Z, P2], then many 
people elaborated on this. For instance, nice updatings and
improvements 
as well as constructions of exotic discrete groups with property (T)
are available in Ballmann-Swiatkowski's paper [BS]. 
That this concerns Kac-Moody groups with compact hyperbolic
Weyl groups was seen independently by Carbone-Garland and Dymara-Januszkiewicz 
[CG, DJ]. According to the first cohomology vanishing criterion [HV,
4.7], the following result [DJ, Corollary 2] proves more than property (T). 

\pec 

{\bf Theorem.---~} \it 
Let $G$ be a Kac-Moody group over $\fq$ with compact hyperbolic Weyl group. Then for $q$
large enough and $1 \leq k \leq n-1$, the continuous cohomology groups 
$H^k_{ct} \bigl( {\rm Aut}(\Delta_\pm),\rho \bigr)$
with coefficients in any unitary representation $\rho$ vanish. 
\cqfd\rm\pec

{\it Remarks.---~} 1. Note that the case for which
we have a positive answer is for a simplex as chamber. 
Actually, negative results are proved for
Fuchsian buildings with non simplicial chambers. 

2. Note also the condition on thickness, inside of which is hidden a spectral 
condition on the Laplace operator. It is actually necessary, in
contrast with the classical situation as we will see in the
forthcoming subsection.  

\pec

{\it Finiteness properties.---~}
A way to define these finiteness properties -- $F_n$ or $FP_n$ -- 
is to use cohomology or actions on complexes, but for the first two degrees it is closely
related to finite generation and finite presentability. 
A criterion defining the class of non-classical Kac-Moody groups for
which positive results are known, involves finiteness of standard
subgroups of the Weyl group. More precisely, we say that the Coxeter
matrix $M$ is {\it $n$-spherical } if $W_J = \langle J \rangle$ is finite for any
subset $J$ with $\#J \leq n$. For instance, if $M$ is of irreducible affine or compact
hyperbolic type, it is $(r-1)$-spherical but not $r$-spherical; and 
2-spherical means that all coefficients are finite.
When results are known, they look like those of arithmetic groups in
characteristic $p$, for which much less is known than in the characteristic 0 case. 
The abstract framework of group actions on twin buildings in [A1]
enables to formulate and prove results in both classical and Kac-Moody cases. 
The arguments there rely on homotopy considerations, in particular make use of a
deep result by K. Brown [Br1]. If we stick to finite generation or
finite presentability, techniques of ordered sets from [T3] can also be used
[AM].  Here are more recent results by P. Abramenko [A2]: \pec

{\bf Theorem.---~} \it 
{\rm (i)} If $\#S > 1$ and $M$ is not $2$-spherical, then $B_\pm$ is
not finitely generated (for arbitrary $q$).

{\rm (ii)} Suppose $q > 3$ and $M$ $2$-spherical but not
$3$-spherical. Then $B_\pm$ is finitely generated but not finitely presentable. 

{\rm (iii)} Suppose that $q > 13$ and $M$ is $3$-spherical.
Then $B_\pm$ is finitely presentable.
\cqfd\rm\pec

{\it Remarks.---~} 1. Non finite generation (resp. non finite presentability) is a way 
to contradict property (T) (resp. Gromov hyperbolicity). 

2. The hypothesis \og $\fq$ large enough\fg once again comes into play. 
P. Abramenko noticed also that over tiny fields, 
some groups as in the theorem are lattices of their building but are 
{\it not } finitely generated. This shows that the condition on thickness is necessary for the
previous theorem concerning property (T), while it is not so in the
classical case (rank argument).

\vfill\eject

\centerline{\bf REFERENCES}

\gec

[A1] P. A{\sevenrm BRAMENKO}: {\it Twin buildings and applications to $S$-arithmetic
groups}.  Springer Lecture Notes in Mathematics {\bf 1641} (1996). 

[A2] P. A{\sevenrm BRAMENKO}: {\it This volume } (1999). 

[AM] P. A{\sevenrm BRAMENKO}, B. M{\sevenrm \"UHLHERR}: 
Pr\'esentations de certaines $BN$-paires jumel\'ees comme sommes amalgam\'ees. 
C. R. Acad. Sci. Paris {\bf 325} (1997), 701-706. 

[Ba] N. B{\sevenrm ARDY}: {\it Syst\`emes de racines infinis}. M\'emoire SMF {\bf 65}
(1996). 

[B$_3$R] V. B{\sevenrm ACK}-V{\sevenrm ALENTE}, 
N. B{\sevenrm ARDY}-P{\sevenrm ANSE}, H. B{\sevenrm EN} M{\sevenrm ESSAOUD}, 
G. R{\sevenrm OUSSEAU}: 
{\it Formes presque d\'eploy\'ees des alg\`ebres de Kac-Moody: classification et racines relatives}. 
J. of Algebra {\bf 171 } (1995), 43-96. 

[Bbk] N. B{\sevenrm OURBAKI}: {\it Groupes et alg\`ebres de Lie \rm IV-VI}. 
Masson (1981). 

[BH] M. B{\sevenrm RIDSON}, A. H{\sevenrm AEFLIGER}: 
{\it Metric spaces of nonpositive curvature}. Springer (1999). 

[Bo] A. B{\sevenrm OREL}: {\it Linear algebraic groups}. Springer GTM {\bf 126} (1990). 

[BoT] A. B{\sevenrm OREL}, J. T{\sevenrm ITS}: {\it Groupes r\'eductifs}. 
Publ. Math. IH\'ES {\bf 27} (1965), 55-151. 

[Bou1] M. B{\sevenrm OURDON}: 
{\it Immeubles hyperboliques, dimension conforme et rigidit\'e de Mostow}. 
GAFA {\bf 7} (1997), 245-268. 

[Bou2] M. B{\sevenrm OURDON}: {\it Sur les immeubles fuchsiens et leur type de quasi-isom\'etrie}. 
Erg. Th. and Dynam. Sys. {\bf 20 } (2000), 343-364. 

[BP1] M. B{\sevenrm OURDON}, H. P{\sevenrm AJOT}: 
{\it Poincar\'e inequalities and quasiconformal structure on the boundary of some hyperbolic buildings}. 
Proc. AMS {\bf 127 } (1999), 2315-2324. 

[BP2] M. B{\sevenrm OURDON}, H. P{\sevenrm AJOT}: 
{\it Rigidity of quasi-isometries for some hyperbolic buildings}.
Commentarii Math. Helv. {\bf 75 } (2000), 701-736. 

[Br1] K. B{\sevenrm ROWN}: {\it Finiteness properties of groups}. 
J. Pure Appl. Algebra {\bf 44 } (1987), 45-75.

[Br2] K. B{\sevenrm ROWN}: {\it Buildings}. Springer (1989). 

[BS] W. B{\sevenrm ALLMANN}, J. S{\sevenrm WIATKOWSKI}: {\it On $L^2$-cohomology and property {\rm (T)} for
automorphism groups of polyhedral cell complexes}. GAFA {\bf 7 } (1997), 615-645. 

[BrT1] F. B{\sevenrm RUHAT}, J. T{\sevenrm ITS}: {\it Groupes r\'eductifs sur un corps local {\rm I. } Donn\'ees
radicielles valu\'ees}. Publ. Math. IH\'ES {\bf 41 } (1972), 5-251. 

[BrT2] F. B{\sevenrm RUHAT}, J. T{\sevenrm ITS}: {\it Groupes r\'eductifs sur un corps local {\rm II. } Sch\'emas
en groupes. Existence d'une donn\'ee radicielle valu\'ee}. 
Publ. Math. IH\'ES {\bf 60 } (1984), 5-184. 

[CG] L. C{\sevenrm ARBONE}, H. G{\sevenrm ARLAND}: {\it Existence and deformations
of lattices in groups of Kac-Moody type}. In preparation. 

[Ch] A. C{\sevenrm HOSSON}: {\it Pr\'esentation des groupes de Kac-Moody, d'apr\`es J.
Tits}. Preprint, University of Picardie -- Amiens (1998). 

[D1] M. D{\sevenrm AVIS}: {\it Groups generated by reflections and aspherical manifolds 
not covered by Euclidean spaces}. Annals of Math. {\bf 117 } (1983), 293-324. 

[D2] M. D{\sevenrm AVIS}: {\it Buildings are {\rm CAT(0)}}. 
pp.108-123 of {\it Geometry and cohomology in group theory}. 
P.H. Kropholler, G.A. Niblo, R. St\"ohr \'editeurs. LMS Lecture Notes
Series {\bf 252 } (1997). 

[DG] M. D{\sevenrm EMAZURE}, P. G{\sevenrm ABRIEL}: {\it Groupes alg\'ebriques}. Masson/North
Holland (1970). 

[DJ] J. D{\sevenrm YMARA}, T. J{\sevenrm ANUSZKIEWICZ}: {\it New Kazhdan groups}. 
Geometri\ae $\,$ Dedicata {\bf 80 } (2000), 311-317. 

[G] H. G{\sevenrm ARLAND}: {\it $p$-adic curvature and the cohomology of discrete
subgroups of $p$-adic groups}. Annals of Math. {\bf 97 } (1973), 375-423. 

[GH] \'E. G{\sevenrm HYS}, P. {\sevenrm DE LA } H{\sevenrm ARPE}: {\it Sur les groupes hyperboliques d'apr\`es Mikhael Gromov}. 
Progress in Mathematics {\bf 83}, Birkh\"auser (1990).  

[GP] D. G{\sevenrm ABORIAU}, F. P{\sevenrm AULIN}: {\it Sur les immeubles hyperboliques}. 
To appear in Geometri\ae $\,$ Dedicata. 

[H1] J.-Y. H{\sevenrm \'EE}: {\it Construction de groupes tordus en 
th\'eorie de Kac-Moody}. C. R. Acad. Sc. Paris {\bf 310 } (1990), 77-80. 

[H2] J.-Y. H{\sevenrm \'EE}: {\it Syst\`emes de racines sur un anneau commutatif 
totalement ordonn\'e}. Geometri\ae $\,$ Dedicata {\bf 37 } (1991), 65-102. 

[H3] J.-Y. H{\sevenrm \'EE}: {\it Sur la torsion de Steinberg-Ree des groupes de
Chevalley et de Kac-Moody}. Th\`ese d'Etat University Paris 11 (1993). 

[Ha] F. H{\sevenrm AGLUND}: {\it Existence, unicit\'e et
homog\'en\'eit\'e de certains immeubles hyperboliques}. Preprint University Paris 11
(1999). 

[HP] F. H{\sevenrm AGLUND}, F. P{\sevenrm AULIN}: {\it Simplicit\'e de groupes d'automorphismes d'espaces \`a courbure n\'egative}. 
(pp.181-248 of {\it The Epstein birthday schrift}. Geom. Topol. Monogr. 1) Geom. Topol.,  
Coventry (1998). 

[Hu1] J. H{\sevenrm UMPHREYS}: {\it Introduction to Lie algebras and representation
theory}.  Springer GTM {\bf 9 } (1972). 

[Hu2] J. H{\sevenrm UMPHREYS}: {\it Reflection groups and Coxeter groups}. 
Cambridge University Press (1990). 

[HV] P. {\sevenrm DE LA } H{\sevenrm ARPE}, A. V{\sevenrm ALETTE}: {\it La propri\'et\'e {\rm (T)} de Kazhdan pour les groupes localement compacts}.
Ast\'erisque {\bf 175 } (1989). 

[K] V. K{\sevenrm AC}: {\it Infinite-dimensional Lie algebras}. 
Cambridge University Press (1990). 

[KL] B. K{\sevenrm LEINER}, B. L{\sevenrm EEB}: {\it Rigidity of quasi-isometries for symmetric spaces and
Euclidian buildings}. Publ. Math. IH\'ES {\bf 86 } (1997), 115-197. 

[KP1] V. K{\sevenrm AC}, D. P{\sevenrm ETERSON}: {\it Regular functions on some infinite-dimensional groups}. 
(pp.141-166 of {\it Arithmetic and geometry}. M. Artin et J. Tate eds). 
Progress in Mathematics {\bf 36}, Birkh\" auser (1983). 

[KP2] V. K{\sevenrm AC}, D. P{\sevenrm ETERSON}: {\it Defining relations for certain infinite-dimensional groups}. 
Ast\'erisque Hors-S\'erie (1984), 165-208. 

[KP3] V. K{\sevenrm AC}, D. P{\sevenrm ETERSON}: {\it On geometric invariant theory for infinite-dimensional groups}. 
Springer Lecture Notes in Mathematics {\bf 1271 } (1987), 109-142. 

[L] \'E. L{\sevenrm EBEAU}: {\it Rigidit\'e et flexibilit\'e de complexes poly\'edraux
\`a courbure n\'egative}. Thesis ENS Lyon (1999). 

[Mar] G.A. M{\sevenrm ARGULIS}: {\it Discrete subgroups of semisimple Lie groups}. 
Springer (1990). 

[Mas] B. M{\sevenrm ASKIT}: {\it Kleinian groups}. Springer (1987). 

[Mat1] O. M{\sevenrm ATHIEU}: {\it Formule des caract\`eres pour les alg\`ebres de
Kac-Moody  g\'en\'erales}. Ast\'erisque {\bf 159-160 } (1988). 

[Mat2] O. M{\sevenrm ATHIEU}: {\it Construction d'un groupe de Kac-Moody et applications}.
Comp. Math. {\bf 69 } (1989), 37-60. 

[Mou] G. M{\sevenrm OUSSONG}: {\it Hyperbolic Coxeter groups}. 
Thesis Ohio State University (1988). 

[MR] B. M{\sevenrm \"UHLHERR},  M. R{\sevenrm ONAN}: {\it Local to global structure in twin buildings}. Invent. Math. {\bf
122 } (1995), 71-81. 

[M\"u] B. M{\sevenrm \"UHLHERR}: {\it Locally split and locally
finite buildings of $2$-spherical type}. J. Reine angew. Math. {\bf
511 } (1999), 119-143. 

[P1] P. P{\sevenrm ANSU}: {\it M\'etriques de Carnot-Carath\'eodory et quasi-isom\'etries
des espaces sym\'etriques de rang $1$}. Annals of Math. {\bf 129 } (1989), 1-60. 

[P2] P. P{\sevenrm ANSU}: {\it Formules de Matsushima, de Garland et propri\'et\'e {\rm
(T) } pour des groupes agissant sur des espaces sym\'etriques ou des immeubles}. 
Bull. Soc. Math. France {\bf 126 } (1998), 107-139. 

[Ra1] J. R{\sevenrm AMAGGE}: {\it A realization of certain affine Kac-Moody groups of types II and III}. 
J. Algebra {\bf 171} (1995), 713-806. 

[Ra2] J. R{\sevenrm AMAGGE}: {\it On certain fixed point subgroups of affine Kac-Moody groups}. 
J. Algebra {\bf 171} (1995), 473-514.

[R\'e1] B. R{\sevenrm \'EMY}: {\it Immeubles \`a courbure n\'egative et th\'eorie de 
Kac-Moody}. Manuscript  available at www.iecn.u-nancy.fr/$^\sim$remy/, Nancy (1998). 

[R\'e2] B. R{\sevenrm \'EMY}: {\it Formes presque d\'eploy\'ees des
groupes de Kac-Moody sur un corps quelconque}. PhD Univ. Nancy (1999), to appear in Ast\'erisque, SMF. 

[R\'e3] B. R{\sevenrm \'EMY}: {\it Construction de r\'eseaux en
th\'eorie de Kac-Moody}. C. R. Acad. Sc. Paris {\bf 329},  475-478 (1999).

[Ron] M. R{\sevenrm ONAN}: {\it Lectures on buildings}. Academic press (1989). 

[Rou1] G. R{\sevenrm OUSSEAU}: {\it Almost split {\bf K}-forms of Kac-Moody algebras}. Infinite-dimensional Lie
algebras and groups, Marseille. Adv. Ser. in Math. Physics {\bf 7}, 
World Scientific publications (1989), 70-85. 

[Rou2] G. R{\sevenrm OUSSEAU}: {\it L'immeuble jumel\'e d'une forme presque d\'eploy\'ee d'une alg\`ebre de
Kac-Moody}. Bull. Soc. Math. Belg. {\bf 42 } (1990), 673-694. 

[Rou3] G. R{\sevenrm OUSSEAU}: {\it On forms of Kac-Moody algebras}. Proc. of Symp. in Pure Math. {\bf 56},
part 2, 393-399 (1994).

[RT1] M. R{\sevenrm ONAN},  J. T{\sevenrm ITS}: {\it Building buildings}. Math. Ann. {\bf 278 } (1987), 291-306.

[RT2] M. R{\sevenrm ONAN},  J. T{\sevenrm ITS}: {\it Twin trees. I}. 
Inv. Math. {\bf 116 } (1994), 463-479.

[RT3] M. R{\sevenrm ONAN},  J. T{\sevenrm ITS}: {\it Twin trees. II. Local structure and a universal
construction}. Israel J. Math. {\bf 109 } (1999), 349-377.

[Se] J.-P. S{\sevenrm ERRE}: {\it Arbres, amalgames, ${\rm SL}_2$}. 
Ast\'erisque {\bf 46 } (1977).

[Sp] T. S{\sevenrm PRINGER}: {\it Linear algebraic groups}. 
Progress in Mathematics {\bf 9}, second edition, Birkh\" auser (1998). 

[T1] J. T{\sevenrm ITS}: {\it Buildings of spherical type and finite $BN$-pairs}. 
Springer Lecture Notes in Mathematics {\bf 386 } (1974). 

[T2] J. T{\sevenrm ITS}: {\it On buildings and their applications}. Proceedings ICM Vancouver (1974), 209-220. 

[T3] J. T{\sevenrm ITS}: {\it Ensembles ordonn\'es, immeubles et sommes amalgam\'ees}. 
Bull. Soc. Math. Belg. {\bf 38 } (1986), 367-387. 

[T4] J. T{\sevenrm ITS}: {\it Uniqueness and presentation of Kac-Moody groups over
fields}.  J. of Alg. {\bf 105} (1987), 542-573. 

[T5] J. T{\sevenrm ITS}: {\it Groupes associ\'es aux alg\`ebres de Kac-Moody}. 
Bourbaki talk $\#{\bf 700}$, Ast\'erisque 177-178 (1989). 

[T6] J. T{\sevenrm ITS}: {\it Th\'eorie des groupes}. 
R\'esum\'e de cours, annuaire du Coll\`ege de France, year 88-89,
pp 81-95 and year 89-90, pp 87-103.

[T7] J. T{\sevenrm ITS}: {\it Twin buildings and groups of Kac-Moody type}. 
LMS Lecture Notes series {\bf 165 } (1992), 249-286. 

[V] E. V{\sevenrm INBERG}: {\it Discrete linear groups generated by reflections}. 
Nauk. USSR, Ser. Mat. {\bf 35 } (1971), 1083-1119. 

[Z] A. Z{\sevenrm UK}: {\it La propri\'et\'e {\rm (T) } de Kazhdan pour les groupes
agissant sur les poly\`edres}. C. R. Acad. Sci. Paris {\bf 323 } (1996),
453-458. 

\gec 

\address

\end